\documentclass[11pt,twoside]{article}
\usepackage{amsmath, amsthm, amscd, amsfonts, amssymb, graphicx, color}

\date{}

\begin{document}

\centerline{}

\centerline {\Large{\bf Fuzzy atomic system and fuzzy $K$-frame in fuzzy Hilbert space}}

%% My definition
\newcommand{\mvec}[1]{\mbox{\bfseries\itshape #1}}
\centerline{}
\centerline{\textbf{Prasenjit Ghosh}}
\centerline{Department of Mathematics,  Barwan N.S. High School(HS),}
\centerline{Barwan, Murshidabad, 742132, West Bengal, India}
\centerline{e-mail: prasenjitpuremath@gmail.com}
\centerline{}
\centerline{\textbf{Jayanta Ghosh}}
\centerline{Department of Mathematics, Manickpur Adarsha Vidyapith,}
\centerline{ Howrah-711309, West Bengal, India}
\centerline{e-mail: ghoshjay$_{-}$04@yahoo.com}
\centerline{}
\centerline{\textbf{T. K. Samanta}}
\centerline{Department of Mathematics, Uluberia College,}
\centerline{Uluberia, Howrah, 711315,  West Bengal, India}
\centerline{e-mail: mumpu$_{-}$tapas5@yahoo.co.in}

\newtheorem{Theorem}{\quad Theorem}[section]

\newtheorem{definition}[Theorem]{\quad Definition}

\newtheorem{theorem}[Theorem]{\quad Theorem}

\newtheorem{remark}[Theorem]{\quad Remark}

\newtheorem{corollary}[Theorem]{\quad Corollary}

\newtheorem{note}[Theorem]{\quad Note}

\newtheorem{lemma}[Theorem]{\quad Lemma}

\newtheorem{example}[Theorem]{\quad Example}

\newtheorem{result}[Theorem]{\quad Result}
\newtheorem{conclusion}[Theorem]{\quad Conclusion}

\newtheorem{proposition}[Theorem]{\quad Proposition}

\begin{abstract}
\textbf{\emph{Atomic system in fuzzy Hilbert space is introduced and the existence of the fuzzy atomic systems for a strongly fuzzy bounded linear operator is studied.\,The notion of a K-frame in fuzzy Hilbert space is presented and some of their characterizations are given.\,We will see that fuzzy frame operator of a fuzzy K-frame in fuzzy Hilbert space is invertible under some sufficient condition and validates this by giving some examples.\,Fuzzy \,$K$-frame property in fuzzy Hilbert space preserve by a strongly fuzzy bounded linear operator is established.\,We will describe stability condition of fuzzy \,$K$-frame in fuzzy Hilbert space under some perturbations.\,We construct new types of fuzzy $K$-frame using fuzzy \,$K$-frame in fuzzy Hilbert space.\,Further, it is seen that scalar combinations and product of two fuzzy \,$K$-frames is also a fuzzy \,$K$-frame in fuzzy Hilbert space.
}}
\end{abstract}
{\bf Keywords:}  \emph{Frame, $K$-frame, fuzzy normed linear space, fuzzy inner product    \\ \smallskip \hspace{2cm} space, fuzzy frame.}

\\

%=====================================
\section{Introduction}
%=====================================
 
In 1984, Katsaras \cite{AKK} introduced the concept of fuzzy norm over a linear space.\,Feblin \cite{FB} gave another idea of fuzzy norm on a linear space in 1992.\,Like this way, several authors have developed the concept of fuzzy norm.\,R. Biswas \cite{RB}, A. M. El-Abyed and H. M. El-Hamouly \cite{AME} were among the first who tried to give a meaningful definition of fuzzy inner product space and associated fuzzy norm function.\,Also those definitions are restricted to the real linear space only.\,Recently, P. Muzamdar and S. K. Samanta \cite{PM} introduced a definition of fuzzy inner product space whose associated fuzzy norm is of Bag and Samanta \cite{TB} type.\,Thereafter, A. Hasankhani, A. Nazari and M. Saheli \cite{AHA} introduced a definition of fuzzy inner product space whose associated fuzzy norm is of Felbin type.

The notion of a frame in Hilbert space was first introduced by Duffin and Schaeffer \cite{Duffin} in connection with some fundamental problem in non-harmonic analysis.\,Thereafter, it was further developed and popularized by Daubechies et al.\,\cite{Daubechies} in 1986.\,A frame is a countable family of elements in a separable Hilbert space which allows for a stable, not necessarily unique, decomposition of an arbitrary element into an expansion of the frame element.\,A sequence \,$\left\{\,f_{\,i}\,\right\}_{i \,=\, 1}^{\infty}$\, in a separable Hilbert space \,$H$\, is called a frame for \,$H$, if there exist positive constants \,$0 \,<\, A \,\leq\, B \,<\, \infty$\, such that
\[ A\; \|\,f\,\|^{\,2} \,\leq\, \sum\limits_{i \,=\, 1}^{\infty}\, \left|\ \left <\,f \,,\, f_{\,i} \, \right >\,\right|^{\,2} \,\leq\, B \,\|\,f\,\|^{\,2}\; \;\text{for all}\; \;f \,\in\, H.\]
The constants \,$A$\, and \,$B$\, are called lower and upper bounds, respectively.\,A frame for a Hilbert space is a generalization of an orthonormal basis and this is such a tool that also allows each vector in this space can be written as a linear combination of elements from the frame but, linear independence among the frame elements is not required.\;Such frames play an important role in Gabor and wavelet analysis.\;Several generalizations of frames  namely, \,$g$-frame, \,$K$-frames etc. have been introduced in recent times.\;$K$-frames for a separable Hilbert spaces were introduced by Lara Gavruta \cite{Gavruta} to study the basic notions about atomic system for a bounded linear operator.\;In recent times, \,$K$-frame was presented to reconstruct elements from the range of a bounded linear operator\,$K$\, in a separable Hilbert space.\;$K$-frames are more generalization than the ordinary frames and many properties of ordinary frames may not holds for such generalization of frames.\,Recently, fuzzy frame in fuzzy Hilbert space was introduced by B. Daraby et al.\cite{BDF}

Throughout this paper,\;$H$\; will denote a separable Hilbert space with the inner product \,$\left <\,\cdot \,,\, \cdot\,\right>$\, and \,$\mathcal{B}\,(\,H\,)$\; denote the space of all bounded linear operator on \,$H$.\;We also denote \,$\mathcal{R}\,(\,T\,)$\; for range set of \,$T$\; where \,$T \,\in\, \mathcal{B}\,(\,H\,)$.

%=====================================
\section{Preliminaries}
%=====================================

\begin{theorem}(\,Douglas' factorization theorem\,)\,{\cite{Douglas}}\label{th1}
Let \,$U,\, V \,\in\, \mathcal{B}\,(\,H\,)$.\;Then the following conditions are equivalent:
\begin{itemize}
\item[$(\,i\,)$]\;\;$\mathcal{R}\,(\,U\,) \,\subseteq\, \mathcal{R}\,(\,V\,)$.
\item[$(\,ii\,)$]\;\;$U\,U^{\,\ast} \,\leq\, \lambda^{\,2}\; V\,V^{\,\ast}$\, for some \,$\lambda \,>\, 0$.
\item[$(\,iii\,)$]\;\;$U \,=\, V\,W$\, for some \,$W \,\in\, \mathcal{B}\,(\,H\,)$.
\end{itemize}
\end{theorem}

\begin{theorem}(\cite{Christensen})\label{thm1.1}
Let \,$H_{\,1},\, H_{\,2}$\; be two Hilbert spaces and \;$U \,:\, H_{\,1} \,\to\, H_{\,2}$\; be a bounded linear operator with closed range \;$\mathcal{R}_{\,U}$.\;Then there exists a bounded linear operator \,$U^{\dagger} \,:\, H_{\,2} \,\to\, H_{\,1}$\, such that \,$U\,U^{\dagger}\,x \,=\, x\; \;\forall\; x \,\in\, \mathcal{R}_{\,U}$.
\end{theorem}

\begin{note}
The operator \,$U^{\dagger}$\, defined in Theorem (\ref{thm1.1}), is called the pseudo-inverse of \,$U$.
\end{note}

\begin{definition}{\cite{Gavruta}}
Let \,$K \,\in\, \mathcal{B}\,(\,H\,)$.\;Then a sequence \,$\{\,f_{\,i}\,\}_{i \,=\, 1}^{\infty}$\, in \,$H$\, is said to be a \,$K$-frame for \,$H$\, if there exist constants \,$0 \,<\, A \,\leq\, B \,<\, \infty$\, such that
\[A \,\left \|\,K^{\,\ast}\,f\, \right \|^{\,2} \,\leq\, \sum\limits^{\infty}_{i \,=\, 1}\, \left |\,\left <\,f \,,\, f_{\,i}\,\right >\,\right|^{\,2} \,\leq\, B\;\left\|\,f\,\right\|^{\,2}\; \;\forall\; f \,\in\, H.\]
\end{definition}

\begin{definition}\cite{TB}
Let \,$U$\, be a linear space over the filed \,$F$.\,A fuzzy subset \,$N$\,of \,$U \,\times\, \mathbb{R}$\, is called a fuzzy norm on \,$U$\, if for all \,$x,\, u \,\in\, U$\, and \,$c \,\in\, F$, the following conditions are satisfied:
\begin{itemize}
\item[$(\,N1\,)$]\,\,$\forall\, t \,\in\, \mathbb{R}$\, with \,$t \,\leq\, 0$, \,$N\,(\,x,\, t\,) \,=\, 0$;
\item[$(\,N2\,)$]\,\,$\forall\, t \,\in\, \mathbb{R}$, \,$t \,>\, 0$, \,$N\,(\,x,\, t\,) \,=\, 1$\, if and only if \,$x \,=\, \theta$;
\item[$(\,N3\,)$]\,\,$\forall\, t \,\in\, \mathbb{R}$, \,$t \,>\, 0$, \,$N\,(\,c\,x,\, t\,) \,=\, N\left(\,x,\,\dfrac{t}{|\,c\,|}\,\right)$\, if \,$c \,\neq\, 0$;
\item[$(\,N4\,)$]\,\,$\forall\, s,\, t \,\in\, \mathbb{R}$, \,$x,\, u \,\in\, U$, \,$N\left(\,x \,+\, u,\, s \,+\, t\,\right) \,\geq\, \min\left\{\,N\left(\,x,\, s \,\right),\, N\left(\,u,\, t\,\right)\,\right\}$;
\item[$(\,N5\,)$]\,\,$N\left(\,x,\, \cdot\,\right)$\, is a non-decreasing function on \,$\mathbb{R}$\, and \,$\lim\limits_{t \,\to\, \infty}\,N\,(\,x,\, t\,) \,=\, 1$.
\end{itemize}
The pair \,$\left(\,U,\, N\,\right)$\, is called a fuzzy normed linear space.
\end{definition}

\begin{theorem}\cite{TB}
Let \,$\left(\,U,\, N\,\right)$\, be a fuzzy normed linear space.\,Assume further that,
\begin{itemize}
\item[$(\,N6\,)$]\,\,$\forall\, t \,>\, 0,\, \,N\,(\,x,\, t\,) \,>\, 0\, \,\Rightarrow\, x \,=\, \theta$.
\end{itemize}
\end{theorem}

Define \,$\left\|\,x\,\right\|_{\alpha} \,=\, \bigwedge\left\{\,t \,>\, 0 \,:\, N\left(\,x,\, t\,\right) \,\geq\, \alpha\,\right\},\, \,\alpha \,\in\, (\,0,\, 1\,)$.\,Then \,$\left\{\,\left\|\,\cdot\,\right\|_{\alpha}\,\right\} \,:\, \alpha \,\in\, (\,0,\, 1\,)$\, is an ascending family of norms on \,$U$\, and they are called \,$\alpha$-norms on \,$U$\, corresponding to the fuzzy norm \,$N$\, on \,$U$.

\begin{definition}\cite{TBS}
Let \,$\left(\,U,\, N\,\right)$\, be a fuzzy normed linear space.\,Let \,$\left\{\,x_{n}\,\right\}$\, be a sequence in \,$U$.\,Then \,$\left\{\,x_{n}\,\right\}$\,is said to be convergent if there exists \,$x \,\in\, U$\, such that \,$\lim\limits_{n \,\to\, \infty}\,N\left(\,x_{n} \,-\, x,\, t\,\right) \,=\, 1$, for all \,$t \,>\, 0$. 
\end{definition}

\begin{theorem}\cite{SCC}
Let \,$\left(\,U,\, N\,\right)$\, be a fuzzy normed linear space satisfying \,$(\,N6\,)$\,and \,$\left\{\,x_{n}\,\right\}$\, be a sequence in \,$U$.\,Then \,$\left\{\,x_{n}\,\right\}$\, converges to \,$x$\, if and only if \,$x_{n} \,\to\, x$\, with respect to \,$\left\|\,\cdot\,\right\|_{\alpha}$\, for all \,$\alpha \,\in\, (\,0,\, 1\,)$.  
\end{theorem}

\begin{definition}\cite{TBS}
Let \,$\left(\,U,\, N\,\right)$\, be a fuzzy normed linear space and \,$\alpha \,\in\, (\,0,\, 1\,)$.\,A sequence \,$\left\{\,x_{n}\,\right\}$\, is said to be \,$\alpha$-convergent in \,$U$\, if there exists \,$x \,\in\, U$\, such that \,$\lim\limits_{n \,\to\, \infty}\,N\left(\,x_{n} \,-\, x,\, t\,\right) \,>\, \alpha$, for all \,$t \,>\, 0$.  
\end{definition}

\begin{theorem}\cite{PMS}
Let \,$\left(\,U,\, N\,\right)$\, be a fuzzy normed linear space satisfying \,$(\,N6\,)$.\,If a sequence \,$\left\{\,x_{n}\,\right\}$\,is an \,$\alpha$-convergent in \,$\left(\,U,\, N\,\right)$, then \,$\left\|\,x_{n} \,-\, x\,\right\|_{\alpha} \,\to\, 0$\, as \,$n \,\to\, \infty$.
\end{theorem}

\begin{definition}\cite{PM}
Let \,$U$\, be a linear space over the field \,$\mathbb{C}$\, of complex numbers.\,Let \,$\mu \,:\, U \,\times\, U \,\times\, \mathbb{C} \,\to\, I \,=\, [\,0,\,1\,]$\, be a mapping such that the following conditions and statements hold:
\begin{itemize}
\item[(\,FIP1\,)]\, for \,$s,\, t \,\in\, \mathbb{C}$, \,$\mu\left(\,x \,+\, y,\, z,\, |\,t\,| \,+\, |\,s\,|\,\right) \,\geq\, \min\left\{\,\mu\left(\,x,\, z,\, |\,t\,|\,\right),\, \mu\left(\,y,\, z,\, |\,s\,|\,\right)\,\right\}$;
\item[(\,FIP2\,)]\, for \,$s,\, t \,\in\, \mathbb{C}$, \,$\mu\left(\,x,\, y,\, |\,s\,t\,|\,\right) \,\geq\, \min\left\{\,\mu\left(\,x,\, x,\, |\,s\,|^{2}\,\right),\, \mu\left(\,y,\, y,\, |\,t\,|^{2}\,\right)\,\right\}$;
\item[(\,FIP3\,)]\, for \,$t \,\in\, \mathbb{C}$, \,$\mu\left(\,x,\, y,\, t\,\right) \,=\, \mu\left(\,y,\, x,\, \overline{t}\,\right)$;
\item[(\,FIP4\,)]\, for \,$\alpha\,(\,\neq\, 0\,),\, t \,\in\, \mathbb{C}$, \,$\mu\left(\,\alpha\,x,\, y,\, t\,\right) \,=\, \mu\left(\,x,\, y,\,  \dfrac{t}{|\,\alpha\,|}\,\right)$;
\item[(\,FIP5\,)]\,\,$\mu\left(\,x,\, x,\, t\,\right) \,=\, 0$,\, for all \,$t \,\in\, \mathbb{C} \,/\, \mathbb{R}^{+}$; 
\item[(\,FIP6\,)]\,\,$\mu\left(\,x,\, x,\, t\,\right) \,=\, 1$, for all \,$t \,>\, 0$\, if and only if \,$x \,=\, \theta$;
\item[(\,FIP7\,)]\,\,$\mu\left(\,x,\, x,\, \cdot\,\right) \,:\, \mathbb{R} \,\to\, I$\, is a monotonic non-decreasing function on \,$\mathbb{R}$\, and \,$\lim\limits_{t \,\to\, \infty}\,\mu\left(\,\alpha\,x,\, x,\, t\,\right) \,=\, 1$.
\end{itemize}
Then \,$\mu$\, is called fuzzy inner product on \,$U$\, and \,$\left(\,U,\, \mu\,\right)$\, is called fuzzy inner product space (\,FIP\,) space.
\end{definition}

\begin{theorem}\cite{PM}
Let \,$\left(\,U,\, \mu\,\right)$\, be a FIP.\,Then
\[N\left(\,x,\, t\,\right)\,=\, \begin{cases}
\mu\left(\,x,\, x,\, t^{2}\,\right) & \text{if}\; \,t \,\in\, \mathbb{R},\, \,t \,>\, 0 \\ 0 & \text{if}\; \,t \,\leq\, 0. \end{cases}\]
is a fuzzy norm on \,$U$.
\end{theorem}
Now if \,$\mu$\, satisfies the following conditions:
\begin{itemize}
\item[(\,FIP8\,)]\,\,$\mu\left(\,x,\, x,\, t^{2}\,\right) \,>\, 0,\, \,\forall\, t \,>\, 0 \,\Rightarrow\, x\,=\, \theta$\, and 
\item[(\,FIP9\,)]\,\,for all \,$x,\, y \,\in\, U$\, and \,$p,\, q \,\in\, \mathbb{R}$,
\[\mu\left(x + y,\, x + y,\, 2q^{2}\right) \bigwedge \mu\left(x - y,\, x - y,\, 2p^{2}\right) \geq \mu\left(x,\, x,\, p^{2}\right)\bigwedge\mu\left(y,\, y,\, q^{2}\right),\]
\end{itemize}
then \,$\left\|\,x\,\right\|_{\alpha} \,=\, \bigwedge\left\{\,t \,>\, 0 \,:\, N\left(\,x,\, t\,\right) \,\geq\, \alpha\,\right\},\, \,\alpha \,\in\, (\,0,\, 1\,)$\, is an ordinary norm satisfying the parallelogram law.\,By using polarization identity, we can get ordinary inner product, called the \,$\left<\,\cdot,\, \cdot\,\right>_{\alpha}$-inner product, as follows:
\[\left<\,x,\, y\,\right>_{\alpha} \,=\, \dfrac{1}{4}\left(\,\left\|\,x \,+\, y\,\right\|^{2}_{\alpha} \,-\, \left\|\,x \,-\, y\,\right\|^{2}_{\alpha}\,\right) \,+\, \dfrac{i}{4}\left(\,\left\|\,x \,+\, i\,y\,\right\|^{2}_{\alpha} \,-\, \left\|\,x \,-\, i\,y\,\right\|^{2}_{\alpha}\,\right),\]for all \,$\alpha \,\in\, (\,0,\, 1\,)$. 

\begin{example}\cite{BDF} 
Let \,$\left(\,U,\, \mu\,\right)$\, be a real inner product space.\,Then the function \,$\mu \,:\, U \,\times\, U \,\times\, \mathbb{C} \,\to\, I \,=\, [\,0,\,1\,]$\, defined by 
\[\mu\left(\,x,\, y,\, t\,\right)\,=\, \begin{cases}
\dfrac{|\,t\,|}{|\,t\,| \,+\, \|\,x\,\|\,\|\,y\,\|} & \text{if}\; \,t \,>\,  \|\,x\,\|\,\|\,y\,\|,\\ 0 & \text{if}\; t \,\leq\,  \|\,x\,\|\,\|\,y\,\|,\\0 & \text{if}\; t \,\in\, \mathbb{C} \,\setminus\, \mathbb{R}^{+}
. \end{cases}\]
is a fuzzy inner product.\,Therefore, every classic inner product induces the fuzzy inner product space.\,Here, \,$\mu$\, satisfies (\,FIP8\,) and (\,FIP9\,).\,Also, it follows that
\[\left\|\,x \,+\, y\,\right\|^{2}_{\alpha} \,+\, \left\|\,x \,-\, y\,\right\|^{2}_{\alpha} \,=\, \dfrac{\alpha}{1 \,-\, \alpha}\left(\,2\,\|\,x\,\|^{2} \,+\, 2\,\|\,y\,\|^{2}\,\right),\, \left<\,x,\, y\,\right>_{\alpha} \,=\, \dfrac{\alpha}{1 \,-\, \alpha}\left<\,x,\, y\,\right>\] and \,$\|\,x\,\|_{\alpha} \,=\, \sqrt{\dfrac{\alpha}{1 - \alpha}}\,\|\,x\,\|$. 
\end{example}

Define \,$l^{2}(\mathbb{N}) \,=\, \left\{\,\left\{\,\beta_{i}\,\right\}_{i = 1}^{\infty} \,:\, \sum\limits_{i \,=\, 1}^{\,\infty}|\,\beta_{i}\,|^{2} \,<\, \infty\,\right\}$\, and the consider fuzzy inner product
\[\mu\left(\,x,\, y,\, t\,\right)\,=\, \begin{cases}
1 & \text{if}\; \,t \,>\,  \|\,x\,\|_{l^{2}(\mathbb{N})}\,\|\,y\,\|_{l^{2}(\mathbb{N})},\\ 0 & \text{if}\; t \,\leq\,  \|\,x\,\|_{l^{2}(\mathbb{N})}\,\|\,y\,\|_{l^{2}(\mathbb{N})}
. \end{cases}\]
Now, it can be easily verified that \,$\left\|\,x\,\right\|_{\alpha} \,=\, \|\,x\,\|_{l^{2}(\mathbb{N})}$, for all \,$\alpha \,\in\, (\,0,\,1\,)$
\begin{definition}\cite{PM}
Let \,$\left(\,U,\, \mu\,\right)$\, be a FIP satisfying (\,FIP8\,).\,The linear space \,$U$\, is said to be level complete if for any \,$\alpha \,\in\, (\,0,\, 1\,)$, every Cauchy sequence converges with respect to \,$\left\|\,\cdot\,\right\|_{\alpha}$\, and it is said to be fuzzy Hilbert space it is level complete. 
\end{definition}

\begin{definition}\cite{TBS}
Let \,$\left(\,U,\, N_{1}\,\right)$\, and \,$\left(\,V,\, N_{2}\,\right)$\, be fuzzy normed linear spaces.\,The mapping \,$T \,:\, U \,\to\, V$\, is said to be strongly fuzzy bounded if and only if there exists a positive real number \,$M$\, such that
\[N_{2}\left(\,T(x),\, s\,\right) \,\geq\, N_{1}\left(\,x,\, \dfrac{s}{M}\,\right),\, \,\forall\, x \,\in\, U,\, \forall\, s \,\in\, \mathbb{R}.\] 
\end{definition}

\begin{definition}\cite{TBS}
Let \,$\left(\,U,\, N_{1}\,\right)$\, and \,$\left(\,V,\, N_{2}\,\right)$\, be fuzzy normed linear spaces.\,The mapping \,$T \,:\, U \,\to\, V$\, is said to be uniformly bounded if there exists a positive real number \,$M$\, such that \,$\left\|\,T\,x\,\right\|^{2}_{\alpha} \,\leq\, M\,\left\|\,x\,\right\|^{1}_{\alpha}$\, for all \,$\alpha \,\in\, (\,0,\, 1\,)$, where  \,$\left\|\,\cdot\,\right\|^{1}_{\alpha}$\, and  \,$\left\|\,\cdot\,\right\|^{2}_{\alpha}$\, are \,$\alpha$-norms on \,$N_{1}$\, and \,$N_{2}$, respectively.  
\end{definition}

\,$\mathcal{B}\,(\,U,\, V\,)$\, denotes the set of all strongly fuzzy bounded linear operators from a fuzzy normed linear spaces \,$\left(\,U,\, N_{1}\,\right)$\, to \,$\left(\,V,\, N_{2}\,\right)$.

\begin{theorem}\cite{TBS}
Let \,$\left(\,U,\, N_{1}\,\right)$\, and \,$\left(\,V,\, N_{2}\,\right)$\, be fuzzy normed linear spaces satisfying \,$(\,N_{6}\,)$.\,The mapping \,$T \,:\, U \,\to\, V$\, is strongly fuzzy bounded if and only if it is uniformly bounded with respect to \,$\alpha$-norms on \,$N_{1}$\, and \,$N_{2}$.  
\end{theorem} 

\begin{definition}\cite{TBS}
Let \,$\left(\,U,\, N_{1}\,\right)$\, and \,$\left(\,V,\, N_{2}\,\right)$\, be fuzzy normed linear spaces satisfying \,$(\,N_{6}\,)$.\,For \,$T \,\in\, \mathcal{B}\,(\,U,\, V\,)$, let 
\[\,\left\|\,T\,\right\|^{'}_{\beta} \,=\, \bigvee\limits_{x \in U, x \neq \theta}\dfrac{\left\|\,T\,x\,\right\|^{2}_{\beta}}{\left\|\,x\,\right\|^{1}_{\beta}},\, \beta \in  (\,0,\, 1\,),\]
and \,$\left\|\,T\,\right\|_{\alpha} \,=\, \bigvee\limits_{\beta \,\leq \alpha}\left\|\,T\,\right\|^{'}_{\beta},\, \alpha \in  (\,0,\, 1\,)$.\,Then \,$\left\{\,\left\|\,\cdot\,\right\|_{\alpha}\,:\, \alpha \in  (\,0,\, 1\,)\,\right\}$\, is an ascending family of norms in \,$ \mathcal{B}\,(\,U,\, V\,)$.
\end{definition}

\begin{definition}\cite{PM}
Let \,$\left(\,U,\, \mu\,\right)$\, be a FIP satisfying (\,FIP8\,) and (\,FIP9\,).\,Now, if \,$x,\, y \,\in\, U$\,are such that \,$\left<\,x,\, y\,\right>_{\alpha} \,=\, 0,\, \,\alpha \in  (\,0,\, 1\,)$, then we say that \,$x,\, y$\, are \,$\alpha$-fuzzy orthogonal to each other and denote it by \,$x \,\perp_{\alpha}\, y$.\,Let \,$M$\, be a subset of \,$U$\, and \,$x \,\in\, U$.\,Now, if  \,$\left<\,x,\, y\,\right>_{\alpha} \,=\, 0$\, for all \,$y \,\in\, M$, then we say that \,$x$\, is \,$\alpha$-fuzzy orthogonal to \,$M$\, and denote it by \,$x \,\perp_{\alpha}\, M$.\,The set of all \,$\alpha$-fuzzy orthogonal elements to \,$M$\, is called \,$\alpha$-fuzzy orthogonal set.    
\end{definition}

\begin{definition}\cite{PM}
Let \,$\left(\,U,\, \mu\,\right)$\, be a FIP satisfying (\,FIP8\,) and (\,FIP9\,).\,Now, if \,$x,\, y \,\in\, U$\,are such that \,$\left<\,x,\, y\,\right>_{\alpha} \,=\, 0$, for all \,$\alpha \in  (\,0,\, 1\,)$, then we say that \,$x,\, y$\, are fuzzy orthogonal to each other and denote it by \,$x \,\perp\, y$.\,Thus \,$x \,\perp\, y$\, if and only if \,$x \,\perp_{\alpha}\, y$, for all \,$\alpha \in  (\,0,\, 1\,)$.\,The set of all elements fuzzy orthogonal to each other is called fuzzy orthogonal set.    
\end{definition}
\begin{definition}\cite{PMS}
Let \,$\left(\,U,\, \mu\,\right)$\, be a FIP satisfying (\,FIP8\,) and (\,FIP9\,).\,An \,$\alpha$-fuzzy orthogonal set \,$M$\, in \,$U$\, is said to be \,$\alpha$-fuzzy orthonormal if for all \,$x,\, y \,\in\, M$ 
\[\left<\,x,\, y\,\right>_{\alpha}\,=\, \begin{cases}
1 & \text{if}\; \,x \,=\, y \\ 0 & \text{if}\; \,x \,\neq\, y, \end{cases}\] 
where \,$\left<\,\cdot,\, \cdot\,\right>_{\alpha}$\, is an inner product induced by \,$\mu$.  
\end{definition}

\begin{definition}\cite{PMS}
Let \,$\left(\,U,\, \mu\,\right)$\, be a FIP satisfying (\,FIP8\,) and (\,FIP9\,).\,An fuzzy orthogonal set \,$M$\, in \,$U$\, is said to be fuzzy orthonormal if for all \,$x,\, y \,\in\, M$\,and for all \,$\alpha \in  (\,0,\, 1\,)$, we have  
\[\left<\,x,\, y\,\right>_{\alpha}\,=\, \begin{cases}
1 & \text{if}\; \,x \,=\, y \\ 0 & \text{if}\; \,x \,\neq\, y, \end{cases}\] 
where \,$\left<\,\cdot,\, \cdot\,\right>_{\alpha}$\, is an inner product induced by \,$\mu$.  
\end{definition}

\begin{theorem}\cite{PMS}
Let \,$\left(\,U,\, \mu\,\right)$\, be a FIP satisfying (\,FIP8\,) and (\,FIP9\,) and \,$\left\{\,e_{\,k}\,\right\}_{k \,=\, 1}^{\,\infty}$\, be a fuzzy orthonormal sequence in \,$U$.\,Then the following statements are equivalent:
\begin{itemize}
\item[(i)]\,$\left\{\,e_{\,k}\,\right\}_{k \,=\, 1}^{\,\infty}$\, is complete fuzzy orthonormal;
\item[(ii)]if \,$x \,\perp\, e_{k}$\, for \,$k \,=\, 1,\, 2,\, \cdots$, then \,$x \,=\, \theta$;
\item[(iii)]For every \,$x \,\in\, U$, \,$x \,=\, \sum\limits_{i \,=\, 1}^{\,\infty}\,\left <\,x,\, e_{\,k}\,\right>_{\alpha}\,e_{k}$, for all \,$\alpha \in  (\,0,\, 1\,)$,\ and hence \,$\left <\,x,\, e_{\,k}\,\right>_{\alpha} \,=\, \left <\,x,\, e_{\,k}\,\right>_{\beta}$, for all \,$\alpha,\, \beta \in  (\,0,\, 1\,)$, i.\,e., \,$x$\, is independent of \,$\alpha$. 
\item[(iv)]For every \,$x \,\in\, U$, \,$\|\,x\,\|^{2}_{\alpha} \,=\, \sum\limits_{i \,=\, 1}^{\,\infty}\,\left|\,\left <\,x,\, e_{\,k}\,\right>_{\alpha}\,\right|$, for all \,$\alpha \in  (\,0,\, 1\,)$,\ and hence \,$\|\,x\,\|^{2}_{\alpha} \,=\, \|\,x\,\|^{2}_{\beta}$, for all \,$\alpha,\, \beta \in  (\,0,\, 1\,)$. 
\end{itemize}
\end{theorem}

\begin{theorem}\cite{BDF}
Let \,$\left(\,U,\, \mu\,\right)$\, be a fuzzy Hilbert space satisfying \,$\left(\,\textit{FIP}\,8\,\right)$\, and \,$\left(\,\textit{FIP}\,9\,\right)$.\,Suppose \,$T$\, be a strongly fuzzy bounded linear operator on \,$U$.\,Then there exists a strongly fuzzy bounded linear operator \,$T^{\ast}$\, on \,$U$\, such that \,$\left <\,x,\, T\,y\,\right>_{\alpha} = \left <\,T^{\ast} x,\,y\,\right>_{\alpha}$, for all \,$x,\,y \in U$ and \,$\alpha \in  (\,0,\, 1\,)$. 
\end{theorem}

$T^{\ast}$ is called the adjoint of \,$T$ which satisfies the following:
\begin{itemize}
\item[(i)]$\left(\,T^{\ast}\,\right)^{\ast} \,=\, T$;\,\,$(ii)$\,\,$\left(\,T_{1} + T_{2}\,\right)^{\ast} = T_{1}^{\ast} + T_{2}^{\ast}$;\,\,$(iii)$\,\,$\left(\,\lambda\,T\,\right)^{\ast} = \overline{\lambda}\,T^{\ast}$, for all \,$\lambda \in \mathbb{C}$
\end{itemize}

\begin{definition}\cite{BDF}
Let \,$\left(\,U,\, \mu\,\right)$\, be a fuzzy Hilbert space satisfying \,$\left(\,\textit{FIP}\,8\,\right)$\, and \,$\left(\,\textit{FIP}\,9\,\right)$.\,A countable family of elements \,$\left\{\,f_{\,i}\,\right\}_{i \,=\, 1}^{\,\infty}$\, in \,$U$\, is said to be a fuzzy frame for \,$U$\, if there exist constants \,$0 \,<\, A \,\leq\, B \,<\, \infty$\, such that
\begin{align}
&A\,\left \|\,f\,\right \|_{\alpha}^{\,2} \,\leq\, \sum\limits_{i \,=\, 1}^{\infty}\, \left|\,\left <\,f,\,  f_{\,i}\,\right >_{\alpha}\,\right|^{\,2} \,\leq\, B\,\left\|\,f\,\right\|_{\alpha}^{\,2},\label{2.em2.11}
\end{align}
for all \,$f \,\in\, U$\, and \,$\alpha \in  (\,0,\, 1\,)$.\,The constants \,$A$\, and \,$B$\, are called fuzzy frame bounds.\,If \,$A \,=\, B$, then it is called a tight fuzzy frame for \,$U$.\,If \,$A \,=\, B \,=\, 1$, then it is called a Parseval fuzzy frame for \,$U$.\,If the right inequality of (\ref{2.em2.11}) is satisfied then\,$\left\{\,f_{\,i}\,\right\}_{i \,=\, 1}^{\,\infty}$\, is called a fuzzy Bessel sequence in \,$U$.
\end{definition}

Let \,$\left\{\,f_{\,i}\,\right\}_{i \,=\, 1}^{\,\infty}$\, be a fuzzy frame for \,$U$.\,Then the strongly fuzzy bounded linear operator \,$T_{F} \,:\, l^{2}(\mathbb{N}) \,\to\, U$\, defined by \,$T_{F}\,\left\{\,\beta_{\,k}\,\right\} \,=\, \sum\limits_{i \,=\, 1}^{\,\infty}\,\beta_{\,i}\,f_{\,i}$\, is called the fuzzy pre-frame operator or fuzzy synthesis operator and its adjoint operator 
\,$T^{\ast}_{F} \,:\, U \,\to\, l^{2}(\mathbb{N})$\, given by \,$T_{F}^{\ast}\,f \,=\, \left\{\,\left <\,f,\,  f_{\,i}\,\right >_{\alpha}\,\right\}_{i \,=\, 1}^{\,\infty}$\, is called the fuzzy analysis operator.\,The operator \,$S_{F} \,:\, U \,\to\, U$\, defined by \,$S_{F}\,f \,=\, T_{F}\,T_{F}^{\ast}\,f \,=\, \sum\limits_{i \,=\, 1}^{\,\infty}\,\left <\,f,\,  f_{\,i}\,\right >_{\alpha}\,f_{i}$\, is called the fuzzy frame operator.\,It is easy to verify that \,$\left <\,S_{F}\,f,\,  f\,\right >_{\alpha} \,=\, \sum\limits_{i \,=\, 1}^{\infty}\, \left|\,\left <\,f,\,  f_{\,i}\,\right >_{\alpha}\,\right|^{\,2}$, for all \,$f \,\in\, U$.\,The fuzzy frame operator \,$S_{F}$\, is strongly fuzzy bounded, invertible and self-adjoint.

%=====================================
\section{Fuzzy $K$-frame in fuzzy Hilbert space}
%=====================================

In this section, we first give some properties of fuzzy inner product spaces.

\begin{proof}
Now, 
\begin{align*}
\left\|\,x \,+\, i\,x\,\right\|_{\alpha} &= \bigwedge\left\{\,t \,>\, 0 \,:\, N\left(\,x \,+\, i\,x,\, t\,\right) \,\geq\, \alpha\,\right\} \\
& = \bigwedge\left\{\,t \,>\, 0 \,:\, \mu\left(\,x \,+\, i\,x,\,x \,+\, i\,x,\, t^{2}\,\right) \,\geq\, \alpha\,\right\} \\
& = \bigwedge\left\{\,t \,>\, 0 \,:\, \mu\left(\,x,\,x \,+\, i\,x,\, \dfrac{t^{2}}{\sqrt{2}}\,\right) \,\geq\, \alpha\,\right\} \\
& = \bigwedge\left\{\,t \,>\, 0 \,:\, \mu\left(\,x \,+\, i\,x,\, x,\,\dfrac{t^{2}}{\sqrt{2}}\,\right) \,\geq\, \alpha\,\right\} \\
& = \bigwedge\left\{\,t \,>\, 0 \,:\, \mu\left(\,x,\, x,\,\dfrac{t^{2}}{\sqrt{2}\,\sqrt{2}}\,\right) \,\geq\, \alpha\,\right\} \\
& = \bigwedge\left\{\,t \,>\, 0 \,:\, N\left(\,x,\,\dfrac{t}{\sqrt{2}}\,\right) \,\geq\, \alpha\,\right\}.
\end{align*}
Similarly, it is easy to verify that 
\begin{align*}
\left\|\,x \,-\, i\,x\,\right\|_{\alpha} &= \bigwedge\left\{\,t \,>\, 0 \,:\, N\left(\,x \,-\, i\,x,\, t\,\right) \,\geq\, \alpha\,\right\} \\
& = \bigwedge\left\{\,t \,>\, 0 \,:\, N\left(\,x,\,\dfrac{t}{\sqrt{2}}\,\right) \,\geq\, \alpha\,\right\}. 
\end{align*}
Thus, \,$\left\|\,x \,+\, i\,x\,\right\|_{\alpha} = \left\|\,x \,-\, i\,x\,\right\|_{\alpha}$.\,Putting $x = y$, in the definition of $\alpha$-fuzzy inner product, we get \,$\left<\,x,\, x\,\right>_{\alpha}$ 
\[ =\ \dfrac{1}{4}\left(\left\|\,x + x\,\right\|^{2}_{\alpha} - \left\|\,x - x\,\right\|^{2}_{\alpha}\right) + \dfrac{i}{4}\left(\,\left\|\,x + i\,x\,\right\|^{2}_{\alpha} - \left\|\,x - i\,x\,\right\|^{2}_{\alpha}\,\right)= \left\|\,x\,\right\|^{2}_{\alpha}.\]
\end{proof}

Since the $\alpha$-norm \,$\left\|\,\cdot\,\right\|_{\alpha}$ is a classic norm on $U$, next we state two theorem whose proof are very similar to the classic Hilbert space.
 
\begin{theorem}\label{th3.1}
Let \,$\left(\,U,\, \mu\,\right)$\, be a fuzzy Hilbert space satisfying \,$\left(\,\textit{FIP}\,8\,\right)$\, and \,$\left(\,\textit{FIP}\,9\,\right)$.\,Let $M$, $N$ be two strongly fuzzy bounded linear operator on $U$.\,Then the following statement are equivalent:
\begin{itemize}
\item[$(\,i\,)$]\;\;$\mathcal{R}\,(\,M\,) \,\subseteq\, \mathcal{R}\,(\,N\,)$.
\item[$(\,ii\,)$]\;\;$M\,M^{\,\ast} \,\leq\, \lambda^{\,2}\; N\,N^{\,\ast}$\, for some \,$\lambda \,>\, 0$.
\item[$(\,iii\,)$]\;\;$M \,=\, N\,W$\, for some strongly fuzzy bounded linear operator \,$W$ on $U$.
\end{itemize}
\end{theorem}

\begin{theorem}\label{th3.2}
Let \,$\left(\,U,\, \mu\,\right)$\, be a fuzzy Hilbert space satisfying \,$\left(\,\textit{FIP}\,8\,\right)$\, and \,$\left(\,\textit{FIP}\,9\,\right)$. Let $T$ be a strongly fuzzy bounded linear operator on $U$ with closed range $\mathcal{R}\,(\,T\,)$.Then there exists a strongly fuzzy bounded linear operator \,$T^{\dagger}$ on $U$ such that \,$T\,T^{\dagger}\,x \,=\, x\; \;\forall\; x \,\in\, \mathcal{R}(\,T)$.\,The operator \,$T^{\dagger}$\, is called the $\alpha$-pseudo-inverse of \,$T$.
\end{theorem}

Now, we give the definition of a fuzzy atomic system for a strongly fuzzy bounded linear operator.

\begin{definition}\label{3.def3.11}
Let \,$\left(\,U,\, \mu\,\right)$\, be a fuzzy Hilbert space satisfying \,$\left(\,\textit{FIP}\,8\,\right)$\, and \,$\left(\,\textit{FIP}\,9\,\right)$.\,Let \,$K$\, be a strongly fuzzy bounded linear operator on \,$U$\, and \,$\left\{\,f_{\,i}\,\right\}_{i \,=\, 1}^{\,\infty}$\, be a sequence of vectors in \,$U$.\;Then \,$\left\{\,f_{\,i}\,\right\}_{i \,=\, 1}^{\,\infty}$\, is said to be a fuzzy atomic system for \,$K$\, if 
\begin{itemize}
\item[$(i)$]\,the series \,$\sum\limits_{i \,=\, 1}^{\,\infty}\,\beta_{\,i}\,f_{\,i}$\, is \,$\alpha$-convergent for all \,$\{\,\beta_{\,i}\,\}_{i \,=\, 1}^{\,\infty} \,\in\, l^{\,2}\,(\,\mathbb{N}\,)$,
\item[$(ii)$] for any \,$f \,\in\, U$, there exists \,$\{\,\beta_{\,i}\,\}_{i \,=\, 1}^{\,\infty} \,\in\, l^{\,2}\,(\,\mathbb{N}\,)$\, such that \,$K\,(\,f\,) =\, \sum\limits_{i \,=\, 1}^{\,\infty}\,\beta_{\,i}\,f_{\,i}$\, and \,$\left\|\,\{\,\beta_{\,i}\,\}_{i \,=\, 1}^{\,\infty}\,\right\|_{l^{2}} \,\leq\, C\, \left\|\,f\,\right\|_{\alpha},$ for some \,$C \,>\, 0$. 
\end{itemize} 
\end{definition}

In the Definition \ref{3.def3.11}, condition \,$(\,i\,)$\, implies that \,$\left\{\,f_{\,i}\,\right\}_{i \,=\, 1}^{\,\infty}$\, is a fuzzy Bessel sequence in \,$U$.

In the following theorem, we present the existence of the fuzzy atomic systems for a strongly fuzzy bounded linear operator.

\begin{theorem}
A separable fuzzy Hilbert space has a fuzzy atomic system for every strongly fuzzy bounded linear operator.
\end{theorem}

\begin{proof}
Let \,$K \,\in\, \mathcal{B}(\,U\,)$\, and \,$\left\{\,e_{\,i}\,\right\}_{i \,=\, 1}^{\,\infty}$\, be a \,$\alpha$-fuzzy orthonornal basis for \,$U$.\,Then \,$f \,=\, \sum\limits_{i \,=\, 1}^{\,\infty}\,\left <\,f,\, e_{\,i}\,\right>_{\alpha}\,e_{i}$\, and therefore \,$K\,f \,=\, \sum\limits_{i \,=\, 1}^{\,\infty}\,\left <\,f,\, e_{\,i}\,\right>_{\alpha}\,K\,e_{i}$.\,Now, we take \,$f_{i} \,=\, K\,e_{i}$\, and \,$\beta_{i} \,=\, \left <\,f,\, e_{\,i}\,\right>_{\alpha}$, for \,$i \,=\, 1,\, 2,\, 3,\, \cdots$.\,Then
\begin{align*}
\sum\limits_{i \,=\, 1}^{\,\infty}\,\left|\,\left<\,f,\, f_{\,i}\,\right>_{\alpha}\,\right|^{\,2} &\,=\, \sum\limits_{i \,=\, 1}^{\,\infty}\,\left|\,\left<\,f,\, K\,e_{\,i}\,\right>_{\alpha}\,\right|^{\,2} \,=\, \sum\limits_{i \,=\, 1}^{\,\infty}\,\left|\,\left<\,K^{\,\ast}\,f,\, e_{\,i}\,\right>_{\alpha}\,\right|^{\,2} \\
&\,=\, \left\|\,K^{\,\ast}\,f\,\right\|_{\alpha}^{\,2} \,\leq\, \left\|\,K^{\,\ast}\,\right\|^{\,2}\,\left\|\,f\,\right\|_{\alpha}^{\,2}.  
\end{align*}  
Also, 
\[\left\|\,\{\,\beta_{\,i}\,\}_{i \,=\, 1}^{\,\infty}\,\right\|_{l^{2}} \,=\, \sum\limits_{i \,=\, 1}^{\,\infty}\,\left|\,\left<\,f,\, e_{\,i}\,\right>_{\alpha}\,\right|^{\,2} \,=\, \left\|\,f\,\right\|_{\alpha}^{\,2}.\]
\end{proof}
Next theorem provides a characterization of fuzzy atomic system.

\begin{theorem}\label{3.tmn2.20}
Let \,$\left\{\,f_{\,i}\,\right\}^{\infty}_{i \,=\, 1}$\, be a sequence in \,$U$\, and \,$K \,\in\, \mathcal{B}\left(\,U\,\right)$.\,Then the following statements are equivalent:
\begin{itemize}
\item[$(i)$]\,$\left\{\,f_{\,i}\,\right\}^{\infty}_{i \,=\, 1}$\, is a fuzzy atomic system for \,$K$\, in \,$U$.
\item[$(ii)$]There exist \,$A,\, B \,>\, 0$\, such that 
\[A \,\left \|\,K^{\ast}\,f\, \right \|_{\alpha}^{\,2} \,\leq\, \sum\limits^{\infty}_{i \,=\, 1}\, \left |\, \left <\,f,\, f_{\,i} \,\right>_{\,\alpha}\,\right |^{\,2} \,\leq\, B \,\left \|\,f\, \right\|_{\alpha}^{\,2},\]for all \,$f \,\in\, U$.
\end{itemize}
\end{theorem}

\begin{proof}$(i)\,\Rightarrow\,(ii)$
Let \,$\left\{\,f_{\,i}\,\right\}_{i \,=\, 1}^{\,\infty}$\, be a fuzzy atomic system for \,$K$\, in \,$U$.\,Then for \,$g \,\in\, U$, there exists \,$\{\,\beta_{\,i}\,\}_{i \,=\, 1}^{\,\infty} \,\in\, l^{\,2}\,(\,\mathbb{N}\,)$\, such that \,$K\,g =\, \sum\limits_{i \,=\, 1}^{\,\infty}\,\beta_{\,i}\,f_{\,i}$\,, where 
\[\left\|\,\{\,\beta_{\,i}\,\}_{i \,=\, 1}^{\,\infty}\,\right\|_{l^{2}} \,\leq\, C\, \left\|\,g \,\right\|_{\alpha},\] for some \,$C \,>\, 0$.  
Now, for each \,$f \,\in\, U$, we have
\begin{align*}
\left \|\,K^{\ast}\,f\,\right \|_{\alpha}^{\,2} &\,=\, \bigvee_{\left\|\,g \,\right\|_{\alpha} \,=\, 1}\,\left|\,\left<\,K^{\ast}\,f,\, g \,\right>_{\alpha}\,\right|^{\,2} \,=\, \bigvee_{\left\|\,g \,\right\|_{\alpha} \,=\, 1}\,\left|\,\left<\,f,\, K\,g \,\right>_{\alpha}\,\right|^{\,2}\\
&\,=\, \bigvee_{\left\|\,g \,\right\|_{\alpha} \,=\, 1}\,\left|\,\left<\,f,\, \sum\limits_{i \,=\, 1}^{\,\infty}\,\beta_{\,i}\,f_{\,i}\,\right>_{\alpha}\,\right|^{\,2} \,=\, \bigvee_{\left\|\,g \,\right\|_{\alpha} \,=\, 1}\,\left|\,\sum\limits_{i \,=\, 1}^{\,\infty}\,\overline{\,\beta_{i}}\left<\,f,\, f_{\,i}\,\right>_{\alpha}\,\right|^{\,2}\\
& \,\leq\, \bigvee_{\left\|\,g \,\right\|_{\alpha} \,=\, 1}\,\sum\limits_{i \,=\, 1}^{\,\infty}\,\left|\,\beta_{\,i}\,\right|^{\,2}\,\sum\limits_{i \,=\, 1}^{\,\infty}\,\left|\,\left<\,f,\, f_{\,i}\,\right>_{\alpha}\,\right|^{\,2}\\
& \,\leq\, \bigvee_{\left\|\,g \,\right\|_{\alpha} \,=\, 1}\,C^{\,2}\, \left\|\,g \,\right\|_{\alpha}^{\,2}\,\sum\limits_{i \,=\, 1}^{\,\infty}\,\left|\,\left<\,f,\, f_{\,i}\,\right>_{\alpha}\,\right|^{\,2}.
\end{align*}
This implies that for each \,$f \,\in\, U$, we have
\[\dfrac{1}{C^{\,2}}\,\left \|\,K^{\ast}\,f\, \right \|_{\alpha}^{\,2} \,\leq\, \sum\limits_{i \,=\, 1}^{\,\infty}\,\left|\,\left<\,f,\, f_{\,i}\,\right>_{\alpha}\,\right|^{\,2}.\]

$(ii)\,\Rightarrow\,(i)$\,\,\,Suppose \,$(ii)$\, holds.\,Then \,$\left\{\,f_{\,i}\,\right\}_{i \,=\, 1}^{\,\infty}$\, is a fuzzy Bessel sequence in \,$U$.\,From the left inequality in \,$(ii)$, it is easy to verify that \,$A\,K\,K^{\,\ast} \,\leq\, T_{F}\,T_{F}^{\,\ast}$.\,Then by the Theorem \ref{th3.1}, there exists a strongly fuzzy bounded linear operator \,$L \,\in\,  \mathcal{B}\left(\,U,\, l^{\,2}\left(\,\mathbb{N}\,\right)\,\right)$\, such that \,$K \,=\, T_{F}\,L$.\,Now, for \,$f \,\in\, U$, we define \,$L\,f \,=\, \,\left\{\,f_{\,i}\,\right\}_{i \,=\, 1}^{\,\infty}$. Then for each \,$f \,\in\, U$, we have 
\[K\,f \,=\, T_{F}\,\left\{\,f_{\,i}\,\right\}_{i \,=\, 1}^{\,\infty} \,=\, \sum\limits_{i \,=\, 1}^{\,\infty}\,\beta_{\,i}\,f_{\,i}\,,\]
\[\left\|\,\left\{\,f_{\,i}\,\right\}_{i \,=\, 1}^{\,\infty}\,\right\|_{l^{\,2}} \,=\, \left\|\,L\,f\,\right\|_{l^{\,2}} \,\leq\, \|\,L\,\|\,\left \|\,f\, \right \|_{\alpha}.\]Thus, \,$\left\{\,f_{\,i}\,\right\}^{\infty}_{i \,=\, 1}$\, is a fuzzy atomic system for \,$K$\, in \,$U$.\,This completes the proof.
\end{proof}

\begin{corollary}
Let \,$\left\{\,f_{\,i}\,\right\}^{\infty}_{i \,=\, 1}$\, be a fuzzy frame for \,$U$\, with the corresponding frame operator \,$S_{F}$.\,Then \,$\left\{\,f_{\,i}\,\right\}^{\infty}_{i \,=\, 1}$\, is a fuzzy atomic system for \,$S_{F}$\, in \,$U$.
\end{corollary}

\begin{proof}
Let \,$A$\, and \,$B$\, be the fuzzy frame bounds of \,$\left\{\,f_{\,i}\,\right\}^{\infty}_{i \,=\, 1}$.\,Since \,$ \mathcal{R}\left(\,T_{F}\,\right) \,=\,H_{F} \,=\, \mathcal{R}\left(\,S_{F}\,\right)$, by Theorem \ref{th3.1}, there exists \,$\alpha \,>\, 0$\, such that \,$\alpha\,S_{F}\,S_{F}^{\,\ast} \,\leq\, T_{F}\,T_{F}^{\,\ast}$.\,Then for each \,$f \,\in\, H_{F}$, we get
\begin{align*}
&\alpha\,\left \|\,S_{F}^{\ast}\,f \right \|_{\alpha}^{\,2} \,\leq\, \left \|\,T_{F}^{\ast}\,f\, \right \|_{\alpha}^{\,2}
=\,\sum\limits_{i \,=\, 1}^{\,\infty}\,\left|\,\left<\,f,\, f_{\,i}\,\right>_{\alpha}\,\right|^{\,2} \,\leq\, B\,\left \|\,f\, \right \|_{\alpha}^{\,2}.
\end{align*} 
Thus, \,$\left\{\,f_{\,i}\,\right\}_{i \,=\, 1}^{\infty}$\, is a \,fuzzy $S_{F}$-frame for \,$U$\, and hence by Theorem \ref{3.tmn2.20}, it is a fuzzy atomic system a for \,$S_{F}$\, in \,$U$.     
\end{proof}
Now, we give the definition of a fuzzy \,$K$-frame in \,$U$, for some \,$K \,\in\, \mathcal{B}\left(\,U\,\right)$.

\begin{definition}\label{3.def3.21}
Let \,$\left(\,U,\, \mu\,\right)$\, be a fuzzy Hilbert space satisfying \,$\left(\,\textit{FIP}\,8\,\right)$\, and \,$\left(\,\textit{FIP}\,9\,\right)$.\,Let \,$K$\, be a strongly fuzzy bounded linear operator on \,$U$\, and \,$\left\{\,f_{\,i}\,\right\}_{i \,=\, 1}^{\,\infty}$\, be a sequence of vectors in \,$U$.\;Then \,$\left\{\,f_{\,i}\,\right\}_{i \,=\, 1}^{\,\infty}$\, is said to be a fuzzy \,$K$-frame for \,$U$\, if there exist constants \,$0 \,<\, A \,\leq\, B \,<\, \infty$\, such that
\begin{align}
&A\,\left \|\,K^{\,\ast}\,f\,\right \|_{\alpha}^{\,2} \,\leq\, \sum\limits_{i \,=\, 1}^{\infty}\, \left|\,\left <\,f,\,  f_{\,i}\,\right >_{\alpha}\,\right|^{\,2} \,\leq\, B\,\left\|\,f\,\right\|_{\alpha}^{\,2},\label{3.em3.11}
\end{align}
for all \,$f \,\in\, U$\, and \,$\alpha \in  (\,0,\, 1\,)$.\,The constants \,$A$\, and \,$B$\, are called fuzzy \,$K$-frame bounds.\,If \,$A \,=\, B$, then it is called a tight fuzzy $K$-frame for \,$U$.\,If \,$A \,=\, B \,=\, 1$, then it is called a Parseval fuzzy $K$-frame for \,$U$.
\end{definition}

\begin{remark}
\begin{itemize}
\item[$(i)$]In particular, if \,$K \,=\, I$\, in (\ref{3.em3.11}), then \,$\left\{\,f_{\,i}\,\right\}_{i \,=\, 1}^{\infty}$\, is a fuzzy frame for \,$U$.
\item[$(ii)$]Let \,$\left\{\,f_{\,i}\,\right\}_{i \,=\,1}^{\infty}$\, be a tight fuzzy $K$-frame for \,$U$\, with bound \,$A$.\,Then for each \,$f \,\in\, U$, we get
\[ \sum\limits^{\infty}_{i \,=\, 1}\, \left|\,\left <\,f,\; \dfrac{1}{\sqrt{A}}\,f_{\,i}\,\right >_{\alpha}\,\right |^{\,2} \,=\,\left\|\,K^{\,\ast}\,f\,\right\|_{\alpha}^{\,2}.\]
Thus, the family \,$\left\{\,\dfrac{1}{\,\sqrt{A}}\,f_{\,i}\,\right\}^{\infty}_{i \,=\, 1}$\, is a Parseval fuzzy $K$-frame for \,$U$.
\item[$(iii)$]By Theorem \ref{3.tmn2.20}, it follows that the sequence \,$\left\{\,f_{\,i}\,\right\}_{i \,=\, 1}^{\infty}$\, is a fuzzy atomic system in \,$U$\, for \,$K \,\in\, \mathcal{B}\left(\,U\,\right)$\, if and only if it is a fuzzy \,$K$-frame for \,$U$.  
\end{itemize}
\end{remark}

\begin{note}
In general, the fuzzy frame operator of a fuzzy \,$K$-frame for \,$U$\, is not invertible.\,But, if \,$K \,\in\, \mathcal{B}\,(\,U\,)$\, has closed range, then \,$S_{F} \,:\, \mathcal{R}\,(\,K\,) \,\to\, S_{F}\,\left(\,\mathcal{R}\,(\,K\,)\,\right)$\, is an invertible operator.\,For \,$f \,\in\, \mathcal{R}\,(\,K\,)$, we have
\begin{align*}
\left\|\,f\,\right\|_{\alpha}^{\,2} &\,=\, \left\|\,(\,K^{\,\dagger}\,)^{\,\ast}\,K^{\,\ast}\,f\,\right\|_{\alpha}^{\,2}\,\leq\, \left\|\,K^{\,\dagger}\,\right\|^{\,2}\,\left\|\,K^{\,\ast}\,f\,\right\|_{\alpha}^{\,2}.
\end{align*}
Therefore, if \,$\left\{\,f_{\,i}\,\right\}_{i \,=\, 1}^{\infty}$\, is a fuzzy \,$K$-frame for \,$U$\, then inequality (\ref{3.em3.11}), can be written as
\[A\, \left\|\,K^{\,\dagger}\,\right\|^{\,-\, 2}\,\left\|\,f\,\right\|_{\alpha}^{\,2} \,\leq\, \left<\,S_{F}\,f,\, f\,\right>_{\alpha}\,\leq\, B\,\left\|\,f\,\right\|_{\alpha}^{\,2}\]
and furthermore for each \,$f \,\in\, S_{F}\,\left(\,\mathcal{R}\,(\,K\,)\,\right)$, we have
\[B^{\,-\, 1}\left\|\,f\,\right\|_{\alpha}^{\,2} \leq \left<\,S^{\,-\, 1}_{F}\,f,\, f\,\right>_{\alpha} \leq A^{\,-\, 1}\left\|\,K^{\,\dagger}\,\right\|^{\,2}\left\|\,f\,\right\|_{\alpha}^{\,2}.\] 
\end{note}

\begin{example}
Let \,$\left(\,U,\, \left<\,\cdot\,\right>\,\right)$\, be a classic Hilbert space and \,$\left\{\,f_{\,i}\,\right\}_{i \,=\, 1}^{\infty}$\, is a \,$K$-frame for \,$U$\, with bounds \,$A$\, and \,$B$, where \,$K$\, be a uniformly bounded linear operator on \,$U$.\,Then for each \,$f \,\in\, U$, we have
\begin{align*}
&A\,\left \|\,K^{\,\ast}\,f\,\right \|^{\,2} \,\leq\, \sum\limits_{i \,=\, 1}^{\infty}\, \left|\,\left <\,f,\,  f_{\,i}\,\right >\,\right|^{\,2} \,\leq\, B\,\left\|\,f\,\right\|^{\,2}.
\end{align*}
Now, for any \,$\alpha \,\in\, (\,0,\,1\,)$, we get
\begin{align*}
&A\,\frac{\alpha}{1 \,-\, \alpha}\,\left \|\,K^{\,\ast}\,f\,\right \|^{\,2} \,\leq\, \sum\limits_{i \,=\, 1}^{\infty}\, \frac{\alpha}{1 \,-\, \alpha}\left|\,\left <\,f,\,  f_{\,i}\,\right >\,\right|^{\,2} \,\leq\, B\,\,\frac{\alpha}{1 \,-\, \alpha}\left\|\,f\,\right\|^{\,2}.
\end{align*}
This implies that
\begin{align*}
&A\,\left \|\,K^{\,\ast}\,f\,\right \|_{\alpha}^{\,2} \,\leq\, \sum\limits_{i \,=\, 1}^{\infty}\, \left|\,\left <\,f,\,  f_{\,i}\,\right >_{\alpha}\,\right|^{\,2} \,\leq\, B\,\left\|\,f\,\right\|_{\alpha}^{\,2}.
\end{align*}
Thus, \,$\left\{\,f_{\,i}\,\right\}_{i \,=\, 1}^{\infty}$\, is a fuzzy \,$K$-frame for \,$U$\, with bounds \,$A$\, and \,$B$.
\end{example}

\begin{example}\label{exp3.1}
Let \,$U \,=\, \mathbb{C}^{3}$\, and \,$\left\{\,e_{\,1},\, e_{2},\, e_{\,3}\,\right\}$\, be the orthonormal basis for \,$U$.\,Define \,$K \,:\, U \,\to\, U$\, by \,$K\,e_{1} \,=\, e_{1}$, \,$K\,e_{2} \,=\, e_{1} \,-\, e_{2}$\, and \,$K\,e_{3} \,=\, e_{1} \,+\, e_{2}$.\,Then for \,$f \,\in\, U$, we get
\begin{align*}
\left \|\,K^{\,\ast}\,f\,\right \|^{\,2} &\,=\, \sum\limits_{i \,=\, 1}^{3}\, \left|\,\left <\,K^{\,\ast}\,f,\,  e_{\,i}\,\right >\,\right|^{\,2}\\
& \,=\, \left|\,\left <\,f,\,  K\,e_{\,1}\,\right >\,\right|^{\,2}\,+\, \left|\,\left <\,f,\,  K\,e_{\,2}\,\right >\,\right|^{\,2} \,+\, \left|\,\left <\,f,\,  K\,e_{\,3}\,\right >\,\right|^{\,2}\\
&=\,3\,\left|\,\left <\,f,\,  e_{\,1}\,\right >\,\right|^{\,2} \,+\, 2\,\left|\,\left <\,f,\,  e_{\,2}\,\right >\,\right|^{\,2} \,\leq\,3\,\left\|\,f\,\right\|^{\,2}.
\end{align*}
Consider \,$\left\{\,f_{\,i}\,\right\}_{i \,=\, 1}^{3} \,=\, \left\{\,2\,e_{\,1},\, \dfrac{e_{2}}{\sqrt{2}},\,  \dfrac{e_{2}}{\sqrt{2}}\,\right\}$.\,Then for \,$f \,\in\, U$, we have
\begin{align*}
\sum\limits_{i \,=\, 1}^{3}\, \left|\,\left <\,f,\,  f_{\,i}\,\right >\,\right|^{\,2} &\,=\, \left|\,\left <\,f,\,  2\,e_{\,1}\,\right >\,\right|^{\,2} \,+\, \left|\,\left <\,f,\, \dfrac{e_{2}}{\sqrt{2}} \,\right >\,\right|^{\,2}  \,+\, \left|\,\left <\,f,\, \dfrac{e_{2}}{\sqrt{2}} \,\right >\,\right|^{\,2}\\
&=\,4\,\left|\,\left <\,f,\,  e_{\,1}\,\right >\,\right|^{\,2} \,+\,\left|\,\left <\,f,\,  e_{\,2}\,\right >\,\right|^{\,2}. 
\end{align*}
This shows that
\begin{align*}
&\,\left \|\,f\,\right \|^{\,2} \,\leq\, \sum\limits_{i \,=\, 1}^{3}\, \left|\,\left <\,f,\,  f_{\,i}\,\right >\,\right|^{\,2} \,\leq\, 4\,\left\|\,f\,\right\|^{\,2}.\\
&\Rightarrow\,\dfrac{1}{3}\,\left \|\,K^{\,\ast}\,f\,\right \|^{\,2} \,\leq\, \sum\limits_{i \,=\, 1}^{3}\, \left|\,\left <\,f,\,  f_{\,i}\,\right >\,\right|^{\,2} \,\leq\, 4\,\left\|\,f\,\right\|^{\,2}.
\end{align*}
Thus, \,$\left\{\,f_{\,i}\,\right\}_{i \,=\, 1}^{3}$\, is a \,$K$-frame for \,$U$\, with bounds \,$1\,/\, 3$\, and \,$4$.
Now, for any \,$\alpha \,\in\, (\,0,\,1\,)$, we get
\begin{align*}
&\,\frac{\alpha}{3\,(\,1 \,-\, \alpha\,)}\,\left \|\,K^{\,\ast}\,f\,\right \|^{\,2} \,\leq\, \sum\limits_{i \,=\, 1}^{\infty}\,\frac{\alpha}{1 \,-\, \alpha} \left|\,\left <\,f,\,  f_{\,i}\,\right >\,\right|^{\,2} \,\leq\, 4\,\,\frac{\alpha}{1 \,-\, \alpha}\left\|\,f\,\right\|^{\,2}.\\
&\Rightarrow\,\dfrac{1}{3}\,\left \|\,K^{\,\ast}\,f\,\right \|_{\alpha}^{\,2} \,\leq\, \sum\limits_{i \,=\, 1}^{3}\, \left|\,\left <\,f,\,  f_{\,i}\,\right >_{\alpha}\,\right|^{\,2} \,\leq\, 4\,\left\|\,f\,\right\|_{\alpha}^{\,2}.
\end{align*}
Hence, \,$\left\{\,f_{\,i}\,\right\}_{i \,=\, 1}^{3}$\, is a fuzzy \,$K$-frame for \,$U$\, with bounds \,$1\,/\, 3$\, and \,$4$.

The fuzzy frame operator \,$S_{F} \,:\, U \,\to\,U$\, is given by
\begin{align*}
S_{F}\,f& \,=\, \sum\limits_{i \,=\, 1}^{3}\, \left <\,f,\,  f_{\,i}\,\right >_{\alpha}\,f_{i}\\
& \,=\, \left <\,f,\,  2\,e_{\,1}\,\right >_{\alpha}\,2\,e_{1} \,+\, \left <\,f,\,  \dfrac{e_{2}}{\sqrt{2}}\,\right >_{\alpha}\,\dfrac{e_{2}}{\sqrt{2}} \,+\, \left <\,f,\,  \dfrac{e_{2}}{\sqrt{2}}\,\right >_{\alpha}\,\dfrac{e_{2}}{\sqrt{2}} \\
&=\, \frac{\alpha}{1 \,-\, \alpha}\,\left[\,4\,\left <\,f,\,  e_{\,1}\,\right >\,e_{1} \,+\, \left <\,f,\,  e_{\,2}\,\right >\,e_{2}\,\right].
\end{align*}
Here, \,$S_{F}\,e_{1} \,=\, \dfrac{4\,\alpha}{1 \,-\, \alpha}\,e_{1}$,\, \,$S_{F}\,e_{2} \,=\, \dfrac{\alpha}{1 \,-\, \alpha}\,e_{2}$\, and \,$S_{F}\,e_{3} \,=\, 0$.
The matrix associated with the operator \,$S_{F,\,G}$\, is given by
\[
\left[\,S_{F}\,\right] \,=\, 
\begin{pmatrix}
\;\dfrac{4\,\alpha}{1 \,-\, \alpha} & 0 & 0\\
\; 0   & \dfrac{\alpha}{1 \,-\, \alpha} & 0\\
\; 0    & 0  & 0
\end{pmatrix}
.\] 
Here, the matrix \,$\left[\,S_{F}\,\right]$\, is not invertible.\,Thus, the operator \,$S_{F}$\, is well defined and non invertible strongly fuzzy bounded linear operator on \,$\mathbb{C}^{\,3}$.
\end{example}

In the following Theorem, we establish a relationship between fuzzy frame and fuzzy \,$K$-frame for \,$U$.

\begin{theorem}
Let \,$K$\, be strongly fuzzy bounded linear operator on \,$U$.\,Then
\begin{itemize}
\item[$(i)$]Every fuzzy frame for \,$U$\, is a fuzzy \,$K$-frame for \,$U$. 
\item[$(ii)$]If \,$\mathcal{R}\left(\,K\,\right)$\, is closed, every fuzzy \,$K$-frame for \,$U$\, is a fuzzy frame for \,$\mathcal{R}\left(\,K\,\right)$. 
\end{itemize}
\end{theorem} 

\begin{proof}$(i)$
Let \,$\left\{\,f_{\,i}\,\right\}_{i \,=\, 1}^{\infty}$\, be a fuzzy frame for \,$U$\, with bounds \,$A$\, and \,$B$.\,Then for each  \,$f \,\in\, U$, we have
\begin{align*}
&\dfrac{A}{\|\,K\,\|^{\,2}}\,\left\|\,K^{\,\ast}\,f\,\right\|_{\alpha}^{\,2} \,\leq\, A\,\left\|\,f\,\right\|_{\alpha}^{\,2}
\leq\, \sum\limits^{\infty}_{i \,=\, 1}\,\left |\, \left <\,f,\, f_{\,i}\,\right >_{\alpha}\,\right |^{\,2} \,\leq\, B\,\left\|\,f\,\, \right\|_{\alpha}^{\,2}.
\end{align*} 
Thus, \,$\left\{\,f_{\,i}\,\right\}_{i \,=\, 1}^{\infty}$\, is a \,fuzzy $K$-frame for \,$U$\, with bounds \,$A \,/\, \|\,K\,\|^{\,2}$\, and \,$B$.

$(ii)$\,\,Let \,$\left\{\,f_{\,i}\,\right\}_{i \,=\, 1}^{\infty}$\, be a fuzzy \,$K$-frame for \,$U$\, with bounds \,$A$\, and \,$B$.\,Since \,$\mathcal{R}\left(\,K\,\right)$\, is closed, by Theorem \ref{th3.2}, there exists an operator \,$K^{\,\dagger} \,\in\, \mathcal{R}\left(\,U\,\right)$\, such that \,$K\,K^{\,\dagger}\,f \,=\, f$, for all \,$f \,\in\, \mathcal{R}\left(\,K\,\right)$.\,Then for each  \,$f \,\in\, \mathcal{R}\left(\,K\,\right)$, we have 
\begin{align*}
&\dfrac{A}{\left\|\,K^{\,\dagger}\,\right\|^{\,2}}\,\left\|\,f\,\right\|_{\alpha}^{\,2} \,\leq\, A\,\left\|\,K^{\,\ast}\,f\,\right\|_{\alpha}^{\,2} \leq\, \sum\limits^{\infty}_{i \,=\, 1}\,\left |\, \left <\,f,\, f_{\,i}\,\right >_{\alpha}\,\right |^{\,2} \,\leq\, B\,\left\|\,f\,\right\|_{\alpha}^{\,2}.
\end{align*} 
Thus, \,$\left\{\,f_{\,i}\,\right\}_{i \,=\, 1}^{\infty}$\, is a  fuzzy frame for \,$\mathcal{R}\left(\,K\,\right)$\, with bounds \,$A \,/\, \left\|\,K^{\,\dagger}\,\right\|^{\,2}$\, and \,$B$. 
\end{proof}

Next, we state a characterization of a fuzzy \,$K$-frame for \,$U$\, with respect to its fuzzy frame operator. 

\begin{remark}\label{3.rmk2.19}
Let \,$\left\{\,f_{\,i}\,\right\}_{i \,=\, 1}^{\infty}$\, be a fuzzy Bessel sequence in \,$U$\, and \,$K \,\in\, \mathcal{B}\left(\,U\,\right)$. Then \,$\left\{\,f_{\,i}\,\right\}_{i \,=\, 1}^{\infty}$\, is fuzzy \,$K$-frame for \,$U$\, if and only if there exists \,$\alpha \,>\, 0$\, such that \,$S_{F} \,\geq\, \alpha\,K\,K^{\,\ast}$, where \,$S_{F}$\, is the corresponding fuzzy frame operator for \,$\left\{\,f_{\,i}\,\right\}_{i \,=\, 1}^{\infty}$.    
\end{remark}

In the next Theorem, we will present some algebric properties of \,fuzzy $K$-frame for \,$U$.  

\begin{theorem}\label{3.tmn2.19}
Let \,$\left\{\,f_{\,i}\,\right\}_{i \,=\, 1}^{\infty}$\, be a \,fuzzy $K_{1}$-frame and fuzzy \,$K_{2}$-frame for \,$U$, where \,$K_{1}$\, and \,$K_{2}$\, be strongly fuzzy bounded linear operators on \,$U$.\,Then for any scalars \,$\alpha,\, \beta$, \,$\left\{\,f_{\,i}\,\right\}_{i \,=\, 1}^{\infty}$\, is fuzzy \,$\alpha\,K_{1} \,+\, \beta\,K_{2}$-frame and fuzzy \,$K_{1}\,K_{2}$-frame for \,$U$.
\end{theorem}

\begin{proof}
Since \,$\left\{\,f_{\,i}\,\right\}_{i \,=\, 1}^{\infty}$\, is a fuzzy \,$K_{1}$-frame and fuzzy \,$K_{2}$-frame for \,$U$, there exist positive constants \,$A_{n},\, B_{n} \,>\, 0$, \,$n \,=\, 1,\, 2$\, such that
\begin{align}
A_{n}\,\left \|\,K_{n}^{\,\ast}\,f\,\right \|_{\alpha}^{\,2}& \,\leq\, \sum\limits_{i \,=\, 1}^{\infty}\, \left|\,\left <\,f,\,  f_{\,i}\,\right >_{\alpha}\,\right|^{\,2} \,\leq\, B_{n}\,\left\|\,f\,\right\|_{\alpha}^{\,2},\label{3.emn2.21}
\end{align}
for all \,$f \,\in\, U$.\,Now, for each \,$f \,\in\, U$, we have
\begin{align*}
&\left \|\,K_{1}^{\,\ast}\,f\,\right \|_{\alpha}^{\,2} \,=\, \dfrac{1}{|\,\alpha\,|^{\,2}}\,\left \|\,\alpha\,K_{1}^{\,\ast}\,f\,\right \|_{\alpha}^{\,2}\\
&=\,\dfrac{1}{|\,\alpha\,|^{\,2}}\,\left \|\,\left(\,\alpha\,K_{1}^{\,\ast} \,+\, \beta\,K_{2}^{\,\ast}\,\right)\,f \,-\, \beta\,K_{2}^{\,\ast}\,f\,\right \|_{\alpha}^{\,2}\\
&\geq\,\dfrac{1}{|\,\alpha\,|^{\,2}}\,\left(\,\left \|\,\left(\,\alpha\,K_{1}^{\,\ast} \,+\, \beta\,K_{2}^{\,\ast}\,\right)\,f\,\right \|_{\alpha}^{\,2} \,-\, \left \|\,\beta\,K_{2}^{\,\ast}\,f\,\right \|_{\alpha}^{\,2}\,\right).
\end{align*}
This implies that
\begin{align*}
\left \|\,\left(\,\alpha\,K_{1}^{\,\ast} \,+\, \beta\,K_{2}^{\,\ast}\,\right)\,f\,\right \|_{\alpha}^{\,2} &\leq\,|\,\alpha\,|^{\,2}\,\left \|\,K_{1}^{\,\ast}\,f\,\right\|_{\alpha}^{\,2} \,+\, |\,\beta\,|^{\,2}\,\left \|\,K_{2}^{\,\ast}\,f\,\right \|_{\alpha}^{\,2}\\
&\leq\,\max\left\{\,|\,\alpha\,|^{\,2},\, |\,\beta\,|^{\,2}\,\right\}\left\{\left \|\,K_{1}^{\,\ast}\,f\,\right \|_{\alpha}^{\,2} \,+\, \left \|\,K_{2}^{\,\ast}\,f\,\right \|_{\alpha}^{\,2}\right\}. 
\end{align*}
Thus, for each \,$f \,\in\, U$, we have
\begin{align*}
&\dfrac{1}{\max\left\{\,|\,\alpha\,|^{\,2},\, |\,\beta\,|^{\,2}\,\right\}}\,\left \|\,\left(\,\alpha\,K_{1}^{\,\ast} \,+\, \beta\,K_{2}^{\,\ast}\,\right)\,f\,\right \|_{\alpha}^{\,2} \\
&\leq\, \left(\,\dfrac{1}{A_{1}} \,+\, \dfrac{1}{A_{2}}\,\right)\,\sum\limits_{i \,=\, 1}^{\infty}\, \left|\,\left <\,f,\,  f_{\,i} \,\right >_{\alpha}\,\right|^{\,2}.
\end{align*} 
On the other hand, from (\ref{3.emn2.21}), we get
\[\sum\limits_{i \,=\, 1}^{\infty}\, \left|\,\left <\,f,\,  f_{\,i} \,\right >_{\alpha}\,\right|^{\,2} \,\leq\, \dfrac{B_{1} \,+\, B_{2}}{2}\,\left \|\,f\,\right \|_{\alpha}^{\,2}\,\, \,\forall\, f \,\in\, U.\]
Hence, \,$\left\{\,f_{\,i}\,\right\}_{i \,=\, 1}^{\infty}$\, is a fuzzy \,$\alpha\,K_{1} \,+\, \beta\,K_{2}$-frame for \,$U$.

Moreover, for each \,$f \,\in\, U$, we have
\begin{align*}
\left \|\,\left(\,K_{1}\,K_{2}\,\right)^{\,\ast}\,f\,\right \|_{\alpha}^{\,2} &\,=\, \left \|\,K_{2}^{\,\ast}\,K_{1}^{\,\ast}\,f\,\right \|_{\alpha}^{\,2}\leq\,\left\|\,K_{2}^{\,\ast}\,\right\|^{\,2}\,\left \|\,K_{1}^{\,\ast}\,f\,\right \|_{\alpha}^{\,2}.
\end{align*} 
Thus, for each \,$f \,\in\, H_{F}$, we have
\begin{align*}
&\dfrac{A_{1}}{\left\|\,K_{2}\,\right\|^{\,2}}\,\left \|\,\left(\,K_{1}\,K_{2}\,\right)^{\,\ast}\,f\,\right\|_{\alpha}^{\,2} \,\leq\, A_{1}\,\left \|\,K_{1}^{\,\ast}\,f\,\right \|_{\alpha}^{\,2}\\
&\leq\, \sum\limits_{i \,=\, 1}^{\infty}\, \left|\,\left <\,f,\,  f_{\,i}\,\right >_{\alpha}\,\right|^{\,2}\, \leq\, B_{1}\,\left \|\,f\,\right \|_{\alpha}^{\,2}.  
\end{align*}
Therefore, \,$\left\{\,f_{\,i}\,\right\}_{i \,=\, 1}^{\infty}$\, is a fuzzy \,$K_{1}\,K_{2}$-frame for \,$U$.\,This completes the proof.      
\end{proof}

Next Theorem is a generalization of the Theorem \ref{3.tmn2.19}. 

\begin{theorem}
Let \,$K_{j} \,\in\, \mathcal{B}\,\left(\,U\,\right), \,\, j \,=\, 1,\, 2,\, \cdots,\, n$\, and \,$\left\{\,f_{\,i}\,\right\}_{i \,=\, 1}^{\infty}$\, be a fuzzy \,$K_{j}$-frame for \,$U$.\,Then
\begin{itemize}
\item[$(i)$]If \,$a_{\,j}$, for \,$j \,=\, 1,\, 2,\, \cdots,\, n$,  are finite collection of scalars, then \,$\left\{\,f_{\,i}\,\right\}_{i \,=\, 1}^{\infty}$\, is a fuzzy \,$\sum\limits^{n}_{j \,=\, 1}\,a_{\,j}\,K_{j}$-frame for \,$U$.
\item[$(ii)$]\,\,$\left\{\,f_{\,i}\,\right\}_{i \,=\, 1}^{\infty}$\, is a fuzzy \,$\displaystyle \prod_{j \,=\, 1}^{\,n}\,K_{j}$-frame for \,$U$. 
\end{itemize}
\end{theorem}

\begin{proof}$(i)$\,\,
Since \,$\left\{\,f_{\,i}\,\right\}_{i \,=\, 1}^{\infty}$\, is a fuzzy $K_{j}$-frame in \,$U$, for all \,$j$, there exist constants \,$A,\,B \,>\, 0$\, such that
\begin{align*}
A\,\left \|\,K_{j}^{\,\ast}\,f\,\right \|_{\alpha}^{\,2} \,\leq\,\sum\limits_{i \,=\, 1}^{\infty}\, \left|\,\left <\,f,\,  f_{\,i}\,\right >_{\alpha}\,\right|^{\,2} \,\leq\, B\,\left\|\,f\,\right\|_{\alpha}^{\,2}.
\end{align*}
Then for each \,$f \,\in\, U$, we have
\begin{align*}
&\dfrac{A}{n\,\max\limits_{j}\,|\,a_{\,j}\,|^{\,2}}\,\left\|\,\left(\,\sum\limits^{n}_{j \,=\, 1}\,a_{\,j}\,K_{\,j}\,\right)^{\,\ast}\,f\,\right\|_{\alpha}^{\,2} \,\leq\, A\,\left\|\,K_{j}^{\,\ast}\,f\,\right\|_{\alpha}^{\,2}\\
&\hspace{.5cm} \,\leq\, \sum\limits_{i \,=\, 1}^{\infty}\, \left|\,\left <\,f,\,  f_{\,i}\,\right >_{\alpha}\,\right|^{\,2} \,\leq\, B\,\left\|\,f\,\right\|_{\alpha}^{\,2}. 
\end{align*}
Thus, \,$\left\{\,f_{\,i}\,\right\}_{i \,=\, 1}^{\infty}$\, is a fuzzy \,$\sum\limits^{n}_{j \,=\, 1}\,a_{\,j}\,K_{j}$-frame for \,$U$.

Proof of \,$(ii)$.\,\, For each \,$f \,\in\, U$, we have
\begin{align*}
&\dfrac{A}{\displaystyle \prod_{j \,=\, 1}^{n \,-\, 1}\,\left\|\,K^{\,\ast}_{j}\,\right\|^{\,2}}\,\left\|\,\left(\,\displaystyle \prod_{j \,=\, 1}^{\,n}\,K_{j}\,\right)^{\,\ast}\,f\,\right\|_{\alpha}^{\,2} \,\leq\, A\,\left\|\,K^{\,\ast}_{n}\,f\,\right\|_{\alpha}^{\,2}\\
&\hspace{.5cm} \,\leq\, \sum\limits_{i \,=\, 1}^{\infty}\, \left|\,\left <\,f,\,  f_{\,i}\,\right >_{\alpha}\,\right|^{\,2} \,\leq\, B\,\left\|\,f\,\right\|_{\alpha}^{\,2}.
\end{align*} 
Thus, \,$\left\{\,f_{\,i}\,\right\}_{i \,=\, 1}^{\infty}$\, is a fuzzy \,$\displaystyle \prod_{j \,=\, 1}^{\,n}\,K_{j}$-frame for \,$U$.\,This completes the proof.  
\end{proof}

\begin{theorem}
Let \,$\left\{\,f_{\,i}\,\right\}^{\infty}_{i \,=\, 1}$\; and \,$\left\{\,g_{\,i}\,\right\}^{\infty}_{i \,=\, 1}$\; be two fuzzy Bessel sequences in \,$U$\, with bounds \,$C$\, and \,$D$, respectively.\;Suppose that \,$T_{F}$\; and \,$T_{F}^{\,\prime}$\; be their fuzzy pre-frame operators such that \,$T_{F}\,\left(\,T_{F}^{\,\prime}\,\right)^{\,\ast} \,=\, K$\,, for some \,$K \,\in\, \mathcal{B}\,(\,U\,)$.\;Then \,$\left\{\,f_{\,i}\,\right\}^{\infty}_{i \,=\, 1}$\, is a fuzzy \,$K$-frame for \,$U$.\,Furthermore, if \,$T_{F}^{\,\prime}\, T_{F}^{\,\ast}\, \,=\, K$\, then \,$\left\{\,g_{\,i}\,\right\}^{\infty}_{i \,=\, 1}$\, is fuzzy \,$K$-frame for \,$U$.
\end{theorem}

\begin{proof}
Since \,$T_{F}\,\left(\,T_{F}^{\,\prime}\,\right)^{\,\ast} \,=\, K$\,, for \,$f \,\in\, U$, we  have 
\[K\,f \,=\, T_{F}\,\left(\,T_{F}^{\,\prime}\,\right)^{\,\ast}\,f \,=\, \sum\limits^{\infty}_{i \,=\, 1}\,\left(\,\left(\,T_{F}^{\,\prime}\,\right)^{\,\ast}\,f\,\right)_{\,i}\,f_{\,i}\,,\]\, where \,$\left(\,\left(\,T_{F}^{\,\prime}\,\right)^{\,\ast}\,f\,\right)_{\,i}$\, is the \,$i$-th coordinate of \,$\left(\,T_{F}^{\,\prime}\,\right)^{\,\ast}f$.\,Then for \,$f \,\in\, U$, we get
\begin{align*}
\left\|\,K^{\,\ast}\,f\,\right\|_{\alpha}^{\,4}& \,=\, \left|\,\left<\,K\,K^{\,\ast}\,f,\, f\,\right>_{\alpha}\,\right|^{\,2}=\, \left|\,\left<\,\sum\limits^{\infty}_{i \,=\, 1}\,\left(\,\left(\,T_{F}^{\,\prime}\,\right)^{\,\ast}\,K^{\,\ast}\,f\,\right)_{\,i}\,f_{\,i}\,\right>_{\alpha}\,\right|^{\,2}\\
&  \,\leq\, \sum\limits^{\infty}_{i \,=\, 1}\,\left|\,\left(\,\left(\,T_{F}^{\,\prime}\,\right)^{\,\ast}\,K^{\,\ast}\,f\,\right)_{\,i}\,\right|^{\,2}\,\sum\limits^{\infty}_{i \,=\, 1}\,|\,\left<\,f,\, f_{\,i}\,\right>_{\alpha}\,|^{\,2}\\
&\leq\, \left\|\, \left(\,T_{F}^{\,\prime}\,\right)^{\,\ast}\,\right\|^{\,2}\,\left\|\,K^{\,\ast}\,f\,\right\|_{\alpha}^{\,2}\,\sum\limits^{\infty}_{i \,=\, 1}\,|\,\left<\,f,\, f_{\,i}\,\right>_{\alpha}\,|^{\,2}\\
&\leq\, D\,\left\|\,K^{\,\ast}\,f\,\right\|_{\alpha}^{\,2}\,\sum\limits^{\infty}_{i \,=\, 1}\,|\,\left<\,f,\, f_{\,i}\,\right>_{\alpha}\,|^{\,2}.\\
&\Rightarrow \dfrac{1}{D}\, \left\|\,K^{\,\ast}\,f\,\right\|_{\alpha}^{\,2} \,\leq\, \sum\limits^{\infty}_{i \,=\, 1}\,|\,\left<\,f,\, f_{\,i}\,\right>_{\alpha}\,|^{\,2}.
\end{align*}
Hence, \,$\left\{\,f_{\,i}\,\right\}^{\infty}_{i \,=\, 1}$\; is a fuzzy \,$K$-frame for \,$U$.\;Similarly, using  \,$T_{F}^{\,\prime}\,T_{F}^{\,\ast} \,=\, K$, it can be shown that \,$\left\{\,g_{\,i}\,\right\}^{\infty}_{i \,=\, 1}$\; is a fuzzy \,$K$-frame for \,$U$\, with the lower bound \,$1 \,/\, C$.\,This completes the proof.     
\end{proof}

\section{Fuzzy $K$-frame and operator}

Let \,$\left\{\,f_{\,i}\,\right\}^{\infty}_{i \,=\, 1}$\; be a fuzzy \,$K$-frame for \,$U$\, with the corresponding frame operator \,$S_{F}$.\,In Example \ref{exp3.1}, we have seen that a fuzzy \,$K$-frame operator is not always invertible.\,In the Next Example, we will see that fuzzy \,$K$-frame operator is invertible and furthermore verify the representation theorem i.\,e., every \,$f \,\in\, U$\, can be represented as \,$f \,=\, \sum\limits_{i \,=\, 1}^{3}\, \left <\,f,\, S_{F}^{\,-\,1}\,f_{\,i}\,\right >_{\alpha}\,f_{i}$\, and \,$f \,=\, \sum\limits_{i \,=\, 1}^{3}\, \left <\,f,\, f_{\,i}\,\right >_{\alpha}\,S_{F}^{\,-\,1}\,f_{i}$.

\begin{example}\label{exp4.1}
Let \,$U \,=\, \mathbb{R}^{3}$\, and \,$\left\{\,e_{\,1},\, e_{2},\, e_{\,3}\,\right\}$\, be the orthonormal basis for \,$U$.\,Define \,$K \,:\, U \,\to\, U$\, by \,$K\,e_{1} \,=\, e_{1}$, \,$K\,e_{2} \,=\, e_{1}$\, and \,$K\,e_{3} \,=\,e_{2}$.\,Then for \,$f \,\in\, U$, we get \,$\left \|\,K^{\,\ast}\,f\,\right \|^{\,2} \,=\, \sum\limits_{i \,=\, 1}^{3}\, \left|\,\left <\,K^{\,\ast}\,f,\,  e_{\,i}\,\right >\,\right|^{\,2} \,\leq\,2\,\left\|\,f\,\right\|^{\,2}$.\,Consider \,$\left\{\,f_{\,i}\,\right\}_{i \,=\, 1}^{3} \,=\, \left\{\,(\,1,\, 1,\, 1\,),\, (\,1,\, -\,1,\, -\,1\,),\,  (\,0,\, 1,\, -\,2\,)\,\right\}$.\,Then for \,$f \,\in\, U$, we have
\begin{align*}
\sum\limits_{i \,=\, 1}^{3}\, \left|\,\left <\,f,\,  f_{\,i}\,\right >\,\right|^{\,2}& \,=\, \left|\,\left <\,(\,x,\, y,\, z\,),\,  (\,1,\, 1,\, 1\,)\,\right >\,\right|^{\,2}  \,+\, \left|\,\left <\,(\,x,\, y,\, z\,),\,  (\,1,\, -\,1,\, -\,1\,)\,\right >\,\right|^{\,2} \\
&+\,\left|\,\left <\,(\,x,\, y,\, z\,),\,  (\,0,\, 1,\, -\,2\,)\,\right >\,\right|^{\,2}\\
&=\,\left(\,x \,+\, y \,+\, z\,\right)^{2} \,+\, \left(\,x \,-\, y \,-\, z\,\right)^{2} \,+\,  \left(\, y \,-\, 2\,z\,\right)^{2}\\
&=\,2\,x^{\,2} \,+\, 3\,y^{\,2} \,+\, 6\,z^{\,2} \,\leq\, 6\,\left(\,x^{\,2} \,+\, y^{\,2} \,+\, z^{\,2}\,\right) \,=\, 6\,\left\|\,f\,\right\|^{\,2}. 
\end{align*}
This shows that
\begin{align*}
&\,2\,\left \|\,f\,\right \|^{\,2} \,\leq\, \sum\limits_{i \,=\, 1}^{3}\, \left|\,\left <\,f,\,  f_{\,i}\,\right >\,\right|^{\,2} \,\leq\, 6\,\left\|\,f\,\right\|^{\,2}.\\
&\Rightarrow\,\left \|\,K^{\,\ast}\,f\,\right \|^{\,2} \,\leq\, \sum\limits_{i \,=\, 1}^{3}\, \left|\,\left <\,f,\,  f_{\,i}\,\right >\,\right|^{\,2} \,\leq\, 6\,\left\|\,f\,\right\|^{\,2}.
\end{align*}
Now, for any \,$\alpha \,\in\, (\,0,\,1\,)$, we get
\begin{align*}
&\,\frac{\alpha}{1 \,-\, \alpha}\,\left \|\,K^{\,\ast}\,f\,\right \|^{\,2} \,\leq\, \sum\limits_{i \,=\, 1}^{\infty}\, \frac{\alpha}{1 \,-\, \alpha}\left|\,\left <\,f,\,  f_{\,i}\,\right >\,\right|^{\,2} \,\leq\, 6\,\frac{\alpha}{1 \,-\, \alpha}\left\|\,f\,\right\|^{\,2}.\\
&\Rightarrow\,\left \|\,K^{\,\ast}\,f\,\right \|_{\alpha}^{\,2} \,\leq\, \sum\limits_{i \,=\, 1}^{3}\, \left|\,\left <\,f,\,  f_{\,i}\,\right >_{\alpha}\,\right|^{\,2} \,\leq\, 6\,\left\|\,f\,\right\|_{\alpha}^{\,2}.
\end{align*}
Hence, \,$\left\{\,f_{\,i}\,\right\}_{i \,=\, 1}^{3}$\, is a fuzzy \,$K$-frame for \,$U$\, with bounds \,$1$\, and \,$6$.

The fuzzy frame operator \,$S_{F} \,:\, U \,\to\,U$\, is given by
\begin{align*}
&S_{F}\,f \,=\, \sum\limits_{i \,=\, 1}^{3}\, \left <\,f,\,  f_{\,i}\,\right >_{\alpha}\,f_{i}\\
&=\,\left <\,(\,x,\, y,\, z\,),\,  (\,1,\, 1,\, 1\,)\,\right >_{\alpha}(\,1,\, 1,\, 1\,)  \,+\, \left <\,(\,x,\, y,\, z\,),\,  (\,1,\, -\,1,\, -\,1\,)\,\right >_{\alpha}(\,1,\, -\,1,\, -\,1\,)\\
&+\,\left <\,(\,x,\, y,\, z\,),\,  (\,0,\, 1,\, -\,2\,)\,\right >_{\alpha}(\,0,\, 1,\, -\,2\,)\\
&=\,\dfrac{\alpha}{1 \,-\, \alpha}\left(\,\left(\,x\,+\,y\,+\,z\,\right)(\,1,\, 1,\, 1\,)\,\right)\,+\,\left(\,\left(\,x\,-\,y\,-\,z\,\right)(\,1,\, -\,1,\, -\,1\,)\,\right)\,+\\
&+\,\dfrac{\alpha}{1 \,-\, \alpha}\,(\,y\,-\,2\,z\,)\,(\,0,\, 1,\, -\,2\,)\\
&=\,\dfrac{\alpha}{1 \,-\, \alpha}\,\left(\,2\,x,\, 3\,y,\, 6\,z\,\right)
\end{align*}
The matrix associated with the operator \,$S_{F,\,G}$\, is given by
\[
\left[\,S_{F}\,\right] \,=\, 
\begin{pmatrix}
\;\dfrac{2\,\alpha}{1 \,-\, \alpha} & 0 & 0\\
\;-\, 0   & \dfrac{3\,\alpha}{1 \,-\, \alpha} & \,-\,0\\
\;-\, 0    & 0  & \dfrac{6\,\alpha}{1 \,-\, \alpha}
\end{pmatrix}
.\] 
Here, \,$\textit{det}\left[\,S_{F}\,\right] \,=\, \dfrac{3\,6\,\alpha^{\,3}}{(\,1 \,-\, \alpha\,)^{\,3}} \,\neq\, 0$.\,Thus, the operator \,$S_{F}$\, is well defined and invertible strongly fuzzy bounded linear operator on \,$\mathbb{R}^{\,3}$.\,Now
\[
\left[\,S_{F}\,\right]^{\,-\, 1} \,=\, \dfrac{(\,1 \,-\, \alpha\,)^{\,3}}{3\,6\,\alpha^{\,3}}
\begin{pmatrix}
\;18\,\left(\dfrac{\alpha}{1 \,-\, \alpha}\,\right)^{\,2} & 0 & 0\\
\;-\, 0   & 12\,\left(\dfrac{\alpha}{1 \,-\, \alpha}\,\right)^{\,2} & \,-\,0\\
\;-\, 0    & 0  & 6\,\left(\dfrac{\alpha}{1 \,-\, \alpha}\,\right)^{\,2}
\end{pmatrix}
.\] 
\[
 \,=\, \dfrac{(\,1 \,-\, \alpha\,)}{3\,6\,\alpha}
\begin{pmatrix}
\;18 & 0 & 0\\
\;-\, 0   & 12\, & \,-\,0\\
\;-\, 0    & 0  & 6\,
\end{pmatrix}
.\] 
Now, for \,$f \,=\, (\,x,\, y,\, z\,) \,\in\, U$, we have
\begin{align*}
&\sum\limits_{i \,=\, 1}^{3}\, \left <\,f,\, S_{F}^{\,-\,1}\,f_{\,i}\,\right >_{\alpha}\,f_{i}\\
&=\,\left <\,(\,x,\, y,\, z\,),\, S_{F}^{\,-\,1}\,(\,1,\, 1,\, 1\,)\,\right >_{\alpha}(\,1,\, 1,\, 1\,)\\
& \,+\, \left <\,(\,x,\, y,\, z\,),\, S_{F}^{\,-\,1}\,(\,1,\, -\,1,\, -\,1\,)\,\right >_{\alpha}(\,1,\, -\,1,\, -\,1\,)\\
&+\,\left <\,(\,x,\, y,\, z\,),\, S_{F}^{\,-\,1}\,(\,0,\, 1,\, -\,2\,)\,\right >_{\alpha}(\,0,\, 1,\, -\,2\,)\\
&=\,\left <\,(\,x,\, y,\, z\,),\, \dfrac{(\,1 \,-\, \alpha\,)}{3\,6\,\alpha}\,(\,18,\, 12,\, 6\,)\,\right >_{\alpha}(\,1,\, 1,\, 1\,)\\
& \,+\, \left <\,(\,x,\, y,\, z\,),\, \dfrac{(\,1 \,-\, \alpha\,)}{3\,6\,\alpha}\,(\,18,\, \,-\,12,\, \,-\,6\,)\,\right >_{\alpha}(\,1,\, -\,1,\, -\,1\,)\\
&+\,\left <\,(\,x,\, y,\, z\,),\, \dfrac{(\,1 \,-\, \alpha\,)}{3\,6\,\alpha}\,(\,0,\, 12,\, \,-\,12\,)\,\right >_{\alpha}(\,0,\, 1,\, -\,2\,)\\
&=\,\dfrac{\alpha}{1 \,-\, \alpha}\left(\,18\,x \,+\, 12\,y \,+\, 6\,z\,\right)\,\dfrac{(\,1 \,-\, \alpha\,)}{3\,6\,\alpha}\,(\,1,\, 1,\, 1\,)\\
&+\,\dfrac{\alpha}{1 \,-\, \alpha}\left(\,18\,x \,-\, 12\,y \,-\, 6\,z\,\right)\,\dfrac{(\,1 \,-\, \alpha\,)}{3\,6\,\alpha}\,(\,1,\, -\,1,\, -\,1\,)\\
&+\,\dfrac{\alpha}{1 \,-\, \alpha}\left(\,12\,y \,-\, 12\,z\,\right)\,\dfrac{(\,1 \,-\, \alpha\,)}{3\,6\,\alpha}\,(\,0,\, 1,\, -\,2\,)\\
&=\,\dfrac{1}{36}\,\left(\,36\,x\,\, 36\,y,\, 36\,z\,\right) \,=\,\left(\,x\,\, y,\, z\,\right) \,=\, f. 
\end{align*}
Similarly, it can be verified that \,$\sum\limits_{i \,=\, 1}^{3}\, \left <\,f,\, f_{\,i}\,\right >_{\alpha}\,S_{F}^{\,-\,1}\,f_{i} \,=\, f,\, \,f \,\in\, U$.\,Thus, the representation theorem is verified in this example.
\end{example}

The following Example illustrate that a fuzzy Bessel sequence \,$\left \{\,f_{\,i}\, \right \}^{\infty}_{i \,=\, 1}$\, is a fuzzy \,$K$-frame but it is not the same for other strongly fuzzy bounded linear operator\,$T$. 

\begin{example}
Let \,$U \,=\, \mathbb{R}^{3}$\, and \,$\left\{\,e_{\,1},\, e_{2},\, e_{\,3}\,\right\}$\, be the orthonormal basis for \,$U$.\,Define \,$K \,:\, U \,\to\, U$\, by \,$K\,e_{1} \,=\, e_{1}$, \,$K\,e_{2} \,=\, e_{1}$\, and \,$K\,e_{3} \,=\,e_{2}$.\,Consider \,$\left\{\,f_{\,i}\,\right\}_{i \,=\, 1}^{3} \,=\, \left\{\,(\,1,\, 1,\, 1\,),\, (\,1,\, -\,1,\, -\,1\,),\,  (\,0,\, 1,\, -\,2\,)\,\right\}$.\,Then by Example \ref{exp4.1}, \,$\left \{\,f_{\,i}\, \right \}^{3}_{i \,=\, 1}$\, is a fuzzy \,$K$-frame for \,$U$\,and corresponding fuzzy frame operator is 
\[
\left[\,S_{F}\,\right] \,=\, 
\begin{pmatrix}
\;\dfrac{2\,\alpha}{1 \,-\, \alpha} & 0 & 0\\
\;-\, 0   & \dfrac{3\,\alpha}{1 \,-\, \alpha} & \,-\,0\\
\;-\, 0    & 0  & \dfrac{6\,\alpha}{1 \,-\, \alpha}
\end{pmatrix}
.\]   
Let \,$T$\, be a strongly fuzzy  bounded linear operator on \,$U$\, which is given by
\[
\left[\,T\,\right] \,=\, 
\begin{pmatrix}
\;\dfrac{2\,\alpha}{1 \,-\, \alpha} & 0 & 0\\
\;-\, 0   & \dfrac{3\,\alpha}{1 \,-\, \alpha} & \,-\,0\\
\;-\, 0    & 1  & \dfrac{6\,\alpha}{1 \,-\, \alpha}
\end{pmatrix}
\] and \,$f \,=\, e_{3} \,\in\, U$.\,Then \,$\sum\limits_{i \,=\, 1}^{3}\, \left|\,\left <\,f,\,  f_{\,i}\,\right >_{\alpha}\,\right|^{\,2} \,=\, 0$.\,On the other hand 
\begin{align*}
\left \|\,T^{\,\ast}\,f\,\right \|_{\alpha}^{\,2}& \,=\, \sum\limits_{i \,=\, 1}^{3}\, \left|\,\left <\,T^{\,\ast}\,f,\, e_{\,i}\,\right >_{\alpha}\,\right|^{\,2} \,=\, \left|\,\left <\,f,\, T\,e_{\,3}\,\right >_{\alpha}\,\right|^{\,2} \\
&=\,\dfrac{\alpha}{1 \,-\, \alpha}\,\left(\,1 \,+\, \dfrac{36\,\alpha^{\,2}}{(\,1 \,-\, \alpha\,)^{\,2}}\,\right). 
\end{align*}
Thus, \,$\left \{\,f_{\,i}\, \right \}^{\infty}_{i \,=\, 1}$\, is not a fuzzy \,$T$-frame for \,$U$.\,Here, note that \,$\mathcal{R}(\,T\,)$\, is not a subspace of \,$\mathcal{R}(\,K\,)$.  
\end{example}

\begin{theorem}\label{thm2}
Let \,$\left \{\,f_{\,i}\, \right \}^{\infty}_{i \,=\, 1}$\, be a fuzzy K-frame for \,$U$\, and \,$T \,\in\, \mathcal{B}\,\left(\,U\,\right)$\, with \,$\mathcal{R}\,(\,T\,) \,\subset\, \mathcal{R}\,(\,K\,)$.\;Then \,$\left \{\,f_{\,i}\, \right \}^{\infty}_{i \,=\, 1}$\; is a fuzzy T-frame for \,$U$.
\end{theorem}

\begin{proof} 
Suppose \,$\left \{\,f_{\,i}\,\right \}^{\infty}_{i \,=\, 1}$\, is a fuzzy \,$K$-frame for \,$U$.\;Then for each \,$f \,\in\, U$, there exist constants \,$A,\, B \,>\, 0$\; such that
\[A\,\left \|\,K^{\,\ast}\,f\,\right \|_{\alpha}^{\,2} \,\leq\, \sum\limits_{i \,=\, 1}^{\infty}\, \left|\,\left <\,f,\,  f_{\,i} \,\right >_{\alpha}\,\right|^{\,2} \,\leq\, B\,\left\|\,f\,\right\|_{\alpha}^{\,2}.\]
Since \,$R\,(\,T\,) \,\subset\, R\,(\,K\,)$, by Theorem \ref{th3.2}, there exists \,$\lambda \,>\, 0$\, such that \,$T\,T^{\,\ast} \,\leq\, \lambda^{\,2}\, K\,K^{\,\ast}$.\;Thus,
\[\dfrac{A}{\lambda^{\,2}}\,\left\|\,T^{\,\ast}\,f\,\right\|_{\alpha}^{\,2} \,=\, \left<\,\dfrac{A}{\lambda^{\,2}}\,T\,T^{\,\ast}\,f,\, f \,\right>_{\alpha} \,\leq\, \left<\,A\;K\,K^{\,\ast}\,f,\, f\,\right>_{\alpha} \,=\, A\,\left \|\,K^{\,\ast}\,f\,\right \|_{\alpha}^{\,2}.\]
Therefore, for each \,$f \,\in\, U$, we get
\[\dfrac{A}{\lambda^{\,2}}\left\|\,T^{\,\ast}\,f\,\right\|_{\alpha}^{\,2} \,\leq\, \sum\limits^{\infty}_{i \,=\, 1} \left |\,\left <\,f,\, f_{\,i}\,\right >_{\alpha}\,\right |^{\,2} \,\leq\, B\,\left\|\,f\,\right\|_{\alpha}^{\,2}.\]
Hence, \,$\left \{\,f_{\,i}\,\right \}^{\infty}_{i \,=\, 1}$\, is a fuzzy \,$T$-frame for \,$U$.
\end{proof}

\begin{theorem}\label{th6}
Let \,$\left \{\,f_{\,i}\, \right \}^{\infty}_{i \,=\, 1}$\, be a fuzzy K-frame for \,$U$\, with bounds \,$A,\,B$\, and \,$T \,\in\, \mathcal{B}\,(\,U\,)$\, be an invertible with \,$T\,K \,=\, K\,T$, then \,$\left \{\,T\,f_{\,i}\,\right \}^{\infty}_{i \,=\, 1}$\, is a fuzzy K-frame for \,$U$.
\end{theorem}

\begin{proof}
Since \,$T$\, is invertible, for each \,$f \,\in\, U$, we get 
\[ \left\|\,K^{\,\ast}\,f\,\right\|_{\alpha}^{\,2} \,=\, \left\|\,\left (\,T^{\,-\, 1}\,\right )^{\,\ast}\,T^{\,\ast}\,K^{\,\ast}\,f\,\right\|_{\alpha}^{\,2} \,\leq\, \left\|\,\left(\,T^{\,-\, 1}\,\right )^{\,\ast}\,\right\|^{\,2} \,\left\|\,T^{\,\ast}\,K^{\,\ast}\,f\,\right \|_{\alpha}^{\,2}.\]
\begin{equation}\label{eq2}
\Rightarrow\; \left \|\,T^{\,-\, 1}\,\right\|^{\,-\, 2}\,\left\|\,K^{\,\ast}\,f\,\right\|_{\alpha}^{\,2} \,\leq\, \left \|\,T^{\,\ast}\,K^{\,\ast}\,f\,\right\|_{\alpha}^{\,2}.
\end{equation}
Now, for \,$f \,\in\, U$, we have  
\[\sum\limits^{\infty}_{i \,=\, 1} \,\left|\,\left <\,f,\, T\,f_{\,i}\,\right >_{\alpha}\,\right |^{\,2} \,=\, \sum\limits^{\infty}_{i \,=\, 1}\,\left|\,\left <\,T^{\,\ast}\,f,\, f_{\,i}\,\right >_{\alpha}\,\right|^{\,2} \,\geq\, A\, \left\|\,K^{\,\ast}\,T^{\,\ast}\,f\,\right\|_{\alpha}^{\,2} \,=\, A\,\left\|\,T^{\,\ast}\,K^{\,\ast}\,f\,\right\|_{\alpha}^{\,2}\]
\[ \geq\, A\;\left \|\,T^{\,-\, 1}\,\right\|^{\,-2}\,\left\|\, K^{\,\ast}\,f\,\right\|_{\alpha}^{\,2}\; \;[\;\text{using} \;(\ref{eq2})].\]
On the other hand, for all \,$f \,\in\, U$, we have
\[\sum\limits^{\infty}_{i \,=\, 1}\,\left|\,\left <\,f,\, T\,f_{\,i}\,\right >_{\alpha}\,\right |^{\,2} \,=\, \sum\limits^{\infty}_{i \,=\, 1}\,\left|\,\left <\,T^{\,\ast}\,f,\, f_{\,i} \,\right >_{\alpha}\,\right|^{\,2} \,\leq\, B\, \left \|\,T^{\,\ast}\,f\,\right\|_{\alpha}^{\,2} \,\leq\, B \,\left\|\,T\, \right \|^{\,2}\;\|\,f\,\|_{\alpha}^{\,2}.\]
Hence, \,$\left \{\,T\,f_{\,i}\,\right \}^{\infty}_{i \,=\, 1}$\; is a fuzzy \,$K$-frame for \,$U$.
\end{proof}

\begin{theorem}
Let \,$\left \{\,f_{\,i}\, \right \}^{\infty}_{i \,=\, 1}$\, be a fuzzy K-frame for \,$U$\, with bounds \,$A,\,B$\, and \,$T \,\in\, \mathcal{B}\,(\,U\,)$\, such that \,$T\,T^{\,\ast} \,=\, I$\, with \,$T \,K \,=\, K \,T$.\;Then \,$\left \{\,T\,f_{\,i}\,\right \}^{\infty}_{i \,=\, 1}$\, is a fuzzy K-frame for \,$U$.
\end{theorem}

\begin{proof} 
Since \,$T\,T^{\,\ast} \,=\, I$, for \,$f \,\in\, U$, we get \,$\left\|\,T^{\,\ast}\,K^{\,\ast}\,f \,\right\|_{\alpha}^{\,2} \,=\, \|\,K^{\,\ast}\,f\,\|_{\alpha}^{\,2}$.\;Since \,$\{\,f_{\,i}\,\}^{\infty}_{i \,=\, 1}$\, is a fuzzy \,$K$-frame, for each \,$f \,\in\, U$, we have 
\[\sum\limits^{\infty}_{i \,=\, 1}\,\left|\,\left <\,f,\, T\, f_{\,i} \,\right >_{\alpha}\,\right|^{\,2} \,=\, \sum\limits^{\infty}_{i \,=\, 1}\,\left|\,\left <\,T^{\,\ast}\,f,\, f_{\,i} \,\right >_{\alpha}\,\right|^{\,2}\hspace{3cm}\]
\[\,\geq\, A\, \left\|\,K^{\,\ast}\,T^{\,\ast}\,f\,\right \|_{\alpha}^{\,2} \,=\, A\, \left\|\,T^{\,\ast}\,K^{\,\ast}\,f \,\right\|_{\alpha}^{\,2} \,=\, A\, \left\|\,K^{\,\ast}\,f\,\right\|_{\alpha}^{\,2}.\]
Thus, we see that \,$\left \{\,T\,f_{\,i}\,\right \}^{\infty}_{i \,=\, 1}$\, satisfies lower fuzzy \,$K$-frame condition.\;Following the proof of the Theorem (\ref{th6}), it can be shown that it also satisfies upper fuzzy \,$K$-frame condition and therefore it is a fuzzy \,$K$-frame for \,$U$.
\end{proof}

\begin{theorem}\label{th7}
Let \,$\left \{\,f_{\,i}\,\right \}^{\infty}_{i \,=\, 1}$\; be a sequence in \,$U$.\;Then \,$\left \{\,f_{\,i}\,\right \}^{\infty}_{i \,=\, 1}$\, is a fuzzy K-frame for \,$U$\, if and only if there exists a strongly fuzzy bounded linear operator \,$T \,:\, l^{\,2}\,(\,\mathbb{N}\,) \,\to\, U$\, such that \,$f_{\,i} \,=\, T\,e_{\,i}$\; and \,${\mathcal{R}}\,(\,K\,) \,\subset\, {\,\mathcal R}\,(\,T\,)$, where \,$\{\,e_{\,i}\,\}^{\infty}_{i \,=\, 1}$\; is an \,$\alpha$-fuzzy orthonormal basis for \,$l^{\,2}\,(\,\mathbb{N}\,)$.
\end{theorem}

\begin{proof}
First we suppose that \,$\{\,f_{\,i}\,\}^{\infty}_{i \,=\, 1}$\, is a fuzzy \,$K$-frame for \,$U$.\,Then, for each \,$f \,\in\, U$, there exist constants \,$A,\, B \,>\, 0$\; such that
\[A\,\left \|\,K^{\,\ast}\,f\,\right \|_{\alpha}^{\,2} \,\leq\, \sum\limits_{i \,=\, 1}^{\infty}\, \left|\,\left <\,f,\,  f_{\,i} \,\right >_{\alpha}\,\right|^{\,2} \,\leq\, B\,\left\|\,f\,\right\|_{\alpha}^{\,2}.\]
Now, we consider the linear operator \,$L \,:\, U \,\to\, l^{\,2}\,(\,\mathbb{N}\,)$\; defined by
\[ L\,(\,f\,) \,=\, \sum\limits^{\infty}_{i \,=\, 1}\,\left <\,f,\, f_{\,i} \,\right >_{\alpha} \,e_{\,i}\; \;\forall\; f \,\in\, U.\]
Since \,$\{\,e_{\,i}\,\}^{\infty}_{i \,=\, 1}$\; is an \,$\alpha$-fuzzy orthonormal basis for \,$l^{\,2}\,(\,\mathbb{N}\,)$, we can write
\[\left\|\,L\,(\,f\,)\,\right\|_{\,l^{\,2}}^{\,2} \,=\, \sum\limits^{\infty}_{i \,=\, 1}\,\left|\,\left <\,f,\, f_{\,i} \,\right >_{\alpha}\,\right|^{\,2} \,\leq\, B\, \|\,f\,\|_{\alpha}^{\,2}.\]
Thus, \,$L$\, is well-defined and strongly fuzzy bounded linear operator on \,$U$.\;So, the adjoint operator \,$L^{\,\ast} \,:\, l^{\,2}\,(\,\mathbb{N}\,) \,\to\, U$\, exists and then for each \,$f \,\in\, U$, we get 
\[\left <\,L^{\,\ast}\,e_{\,i},\, f \,\right >_{\alpha} \,=\, \left <\,e_{\,i},\, L\,(\,f\,)\,\right>_{\alpha} \;=\, \left <\,e_{\,i},\, \sum\limits^{\infty}_{i \,=\, 1}\,\left <\,f,\, f_{\,i}\,\right >_{\alpha} \,e_{\,i}\,\right>_{\alpha}\,=\, \overline{\left <\,f,\, f_{\,i} \,\right>_{\alpha}} \,=\, \left <\,f_{\,i},\, f\,\right>_{\alpha}.\]
The above calculation shows that, \,$L^{\,\ast}\,(\,e_{\,i}\,) \,=\, f_{\,i}$.\;On the other hand \,$A\, \|\,K^{\,\ast}\,f\,\|_{\alpha}^{\,2} \,\leq\, \left\|\,L\,(\,f\,)\,\right \|_{\,l^{\,2}}^{\,2}$\, and this implies that
\[\left<\,A\,K\,K^{\,\ast}\,f,\, f \,\right>_{\alpha} \,\leq\, \left<\,L^{\,\ast}\,L\,f,\, f\,\right>_{\alpha}\Rightarrow\; A \,K\,K^{\,\ast} \,\leq\, T\,T^{\,\ast},\; \text{where}\;\; T \,=\, L^{\,\ast}\] and hence from the  Theorem \ref{th3.2}, \,${\,\mathcal R}\,(\,K\,) \,\subset\, {\,\mathcal R}\,(\,T\,)$.

\text{Conversely}, suppose that \,$T \,:\, l^{\,2}\,(\,\mathbb{N}\,) \,\to\, U$\; be a strongly fuzzy bounded linear operator such that \,$f_{\,i} \,=\, T\,e_{\,i}$\; and \,${\,\mathcal R}\,(\,K\,) \,\subset\, {\,\mathcal R}\,(\,T\,)$.\;We have to show that \,$\{\,f_{\,i}\,\}_{i \,=\, 1}^{\infty}$\; is a fuzzy \,$K$-frame for\,$U$.\;Let \,$g \,\in\, l^{\,2}\,(\,\mathbb{N}\,)$\, then \,$g \,=\, \sum\limits^{\infty}_{i \,=\, 1}\,c_{\,i}\,e_{\,i}$, where \,$c_{\,i} \,=\, \left<\,g,\, e_{\,i} \,\right>_{\alpha}$.\;Now, for all \,$g \,\in\, l^{\,2}\,(\,\mathbb{N}\,)$, we have 
\[\left<\,T^{\,\ast}\,f,\, g \,\right>_{\alpha} \,=\, \left<\,T^{\,\ast}\,f,\, \sum\limits^{\infty}_{i \,=\, 1}\,c_{\,i}\,e_{\,i}\,\right>_{\alpha} \,=\, \sum\limits^{\infty}_{i \,=\, 1}\, \overline{c_{\,i}}\,\left<\,f,\, T\,e_{\,i}\,\right>_{\alpha}\,=\, \sum\limits^{\infty}_{i \,=\, 1}\,\overline{\,c_{\,i}}\,\left<\,f,\, f_{\,i} \,\right>_{\alpha}\] 
\[\,=\, \sum\limits^{\infty}_{i \,=\, 1}\,\overline{\left<\,g,\, e_{\,i} \,\right>_{\alpha}}\,\left<\,f,\, f_{\,i} \,\right>_{\alpha}\,=\, \sum\limits^{\infty}_{i \,=\, 1}\,\left<\,e_{\,i},\, g  \,\right>_{\alpha} \,\left<\,f,\, f_{\,i}\,\right>_{\alpha} \;=\; \left<\,\sum\limits^{\infty}_{i \,=\, 1}\,\left<\,f,\, f_{\,i} \,\right>_{\alpha}\,e_{\,i},\; g  \,\right>_{\alpha}.\] 
This shows that \,$T^{\,\ast}\,(\,f\,) \;=\; \sum\limits^{\infty}_{i \,=\, 1}\,\left<\,f,\, f_{\,i} \,\right>_{\alpha}\,e_{\,i}\; \;\forall\; f \,\in\, U$.\,Thus, for all \,$f \,\in\, U$,
\[\sum\limits^{\infty}_{i \,=\, 1}\,\left|\,\left <\,f,\, f_{\,i} \,\right >_{\alpha}\,\right |^{\,2} \,=\, \sum\limits^{\infty}_{i \,=\, 1}\,\left|\,\left <\,T^{\,\ast}\,f,\, e_{\,i} \,\right >_{\alpha}\,\right |^{\,2} \,=\, \left \|\,T^{\,\ast}\,f\,\right \|_{\alpha}^{\,2}\,\leq\,  \left \|\,T\,\right \|^{\,2}\,\|\,f\,\|_{\alpha}^{\,2}\] 
Thus, \,$\{\,f_{\,i}\,\}^{\infty}_{i \,=\, 1}$\; is a fuzzy Bessel sequence in \,$U$.\;Since \,${\,\mathcal R}\,(\,K\,) \,\subset\, {\,\mathcal R}\,(\,T\,)$, from Theorem \ref{th3.2}, there exists \,$A \,>\, 0$\; such that \,$A\,K\,K^{\,\ast} \,\leq\, T\,T^{\,\ast}$.\;Hence following the proof of the Theorem \ref{thm2}, for all \,$f \,\in\, U$, we get 
\[A\,\|\,K^{\,\ast}\,f\,\|_{\alpha}^{\,2} \,\leq\, \left\|\,T^{\,\ast}\,f\,\right\|_{\alpha}^{\,2} \,=\, \sum\limits^{\infty}_{i \,=\, 1}\,\left | \,\left <\,f,\, f_{\,i} \,\right >_{\alpha}\,\right |^{\,2}.\]
Hence, \,$\{\,f_{\,i}\,\}^{\infty}_{i \,=\, 1}$\; is a fuzzy \,$K$-frame for \,$U$.\,This completes the proof.
\end{proof}

\section{Stability of fuzzy $K$-frame}

Now, we would complete this section by analyzing stability conditions of fuzzy \,$K$-frame for \,$U$\, under some perturbations. 

\begin{theorem}\label{3.tmn2.21}
Let \,$K_{1} \,\in\, \mathcal{B}\,(\,U\,)$\, and \,$\left\{\,f_{\,i}\,\right\}_{i \,=\, 1}^{\infty}$\, be a fuzzy \,$K_{1}$-frame for \,$U$.\,Suppose \,$K_{2} \,\in\, \mathcal{B}\,(\,U\,)$\, and \,$\lambda_{1},\, \lambda_{2} \,\geq\, 0$\, with \,$\lambda_{2} \,<\, 1$\, such that
\begin{align*}
&\left\|\,\left(\,K_{1}^{\,\ast} \,-\, K_{2}^{\,\ast}\,\right)\,f\,\right \|_{\alpha} \,\leq\, \lambda_{1}\,\left \|\,K_{1}^{\,\ast}\,f\,\right \|_{\alpha} \,+\, \lambda_{2}\,\left \|\,K_{2}^{\,\ast}\,f\,\right \|_{\alpha}, 
\end{align*}
for all \,$f \,\in\, U$.\,Then \,$\left\{\,f_{\,i}\,\right\}_{i \,=\, 1}^{\infty}$\, is a fuzzy \,$K_{2}$-frame for \,$U$.  
\end{theorem}

\begin{proof}
Let \,$\left\{\,f_{\,i}\,\right\}_{i \,=\, 1}^{\infty}$\, be a fuzzy \,$K_{1}$-frame for \,$U$\, with bounds \,$A$\, and \,$B$.\,Now, for each \,$f \,\in\, U$, we get
\begin{align*}
&\left \|\,K_{2}^{\,\ast}\,f\,\right \|_{\alpha} \,\leq\, \left \|\,\left(\,K_{1}^{\,\ast} \,-\, K_{2}^{\,\ast}\,\right)\,f\,\right \|_{\alpha} \,+\, \,\left \|\,K_{1}^{\,\ast}\,f\,\right \|_{\alpha} \\
&\leq\, \left(\,1 \,+\, \lambda_{1}\,\right)\,\left \|\,K_{1}^{\,\ast}\,f\,\right \|_{\alpha} \,+\, \lambda_{2}\,\left \|\,K_{2}^{\,\ast}\,f\,\right \|_{\alpha}\\
&\Rightarrow\, \left(\,1 \,-\, \lambda_{2}\,\right)\,\left \|\,K_{2}^{\,\ast}\,f\,\right \|_{\alpha} \,\leq\,  \left(\,1 \,+\, \lambda_{1}\,\right)\,\left \|\,K_{1}^{\,\ast}\,f\,\right\|_{\alpha}.  
\end{align*}
Thus, for each \,$f \,\in\, U$, we have
\begin{align*}
&A\,\left(\,\dfrac{1 \,-\, \lambda_{2}}{1 \,+\, \lambda_{1}}\,\right)^{\,2}\,\left \|\,K_{2}^{\,\ast}\,f\,\right\|_{\alpha}^{\,2} \,\leq\, A\,\left \|\,K_{1}^{\,\ast}\,f\,\right \|_{\alpha}^{\,2}\leq\, \sum\limits_{i \,=\, 1}^{\infty}\, \left|\,\left <\,f,\,  f_{\,i}\,\right >_{\alpha}\,\right|^{\,2}\, \leq\, B\,\left \|\,f\,\right \|_{\alpha}^{\,2}.
\end{align*} 
Hence, \,$\left\{\,f_{\,i}\,\right\}_{i \,=\, 1}^{\infty}$\, is a fuzzy \,$K_{2}$-frame for \,$U$.\,This completes the proof.   
\end{proof}

\begin{corollary}
Let \,$K_{1},\, K_{2} \,\in\, \mathcal{B}\,(\,U\,)$\, and \,$0 \,\leq\, \lambda_{1},\, \lambda_{2} \,<\, 1$\, such that
\begin{align*}
&\left \|\,\left(\,K_{1}^{\,\ast} \,-\, K_{2}^{\,\ast}\,\right)\,f\,\right \|_{\alpha} \,\leq\, \lambda_{1}\,\left \|\,K_{1}^{\,\ast}\,f\,\right \|_{\alpha} \,+\,  \lambda_{2}\,\left \|\,K_{2}^{\,\ast}\,f\,\right \|_{\alpha}, 
\end{align*}
for all \,$f \,\in\, U$.\,Then \,$\left\{\,f_{\,i}\,\right\}_{i \,=\, 1}^{\infty}$\, is a fuzzy \,$K_{1}$-frame for \,$U$\, if and only if it is a fuzzy \,$K_{2}$-frame for \,$U$.
\end{corollary}

\begin{proof}
Proof of this Corollary directly follows from the proof of the Theorem \ref{3.tmn2.21} by interchanging the role of \,$K_{1}$\, and \,$K_{2}$.
\end{proof}

\begin{corollary}
Let \,$K \,\in\, \mathcal{B}\,(\,U\,)$\, and \,$\left\{\,f_{\,i}\,\right\}_{i \,=\, 1}^{\infty}$\, be a fuzzy \,$K$-frame for \,$U$. Suppose \,$0 \,\leq\, \lambda_{1},\, \lambda_{2} \,<\, 1$\, such that
\begin{align*}
&\left \|\,K^{\,\ast}\,f \,-\, f\,\right \|_{\alpha} \,\leq\, \lambda_{1}\,\left \|\,K^{\,\ast}\,f\,\right \|_{\alpha} \,+\, \lambda_{2}\,\left \|\,f\,\right \|_{\alpha}, 
\end{align*}
for all \,$f \,\in\, U$.\,Then \,$\left\{\,f_{\,i}\,\right\}_{i \,=\, 1}^{\infty}$\, is a fuzzy frame for \,$U$.
\end{corollary}

\begin{proof}
Proof of this Corollary directly follows from the Theorem \ref{3.tmn2.21}  by putting \,$K_{1} \,=\, K$\, and \,$K_{2} \,=\, I_{F}$.    
\end{proof}

\begin{theorem}\label{3.tmn2.201}
Let \,$\left\{\,f_{\,i}\,\right\}_{i \,=\, 1}^{\infty}$\, be a fuzzy \,$K$-frame for \,$U$\, with bounds \,$A$\, and \,$B$.\,If there exists a constant \,$M \,>\, 1$\, such that
\begin{align*}
&\sum\limits_{i \,=\, 1}^{\infty}\, \left|\,\left <\,f,\, f_{\,i} \,-\, g_{\,i}\,\right >_{\alpha}\,\right|^{\,2}\leq\, M\, \min\left\{\,\sum\limits_{i \,=\, 1}^{\infty}\, \left|\,\left <\,f,\, f_{\,i}\,\right >_{\alpha}\,\right|^{\,2},\, \sum\limits_{i \,=\, 1}^{\infty}\, \left|\,\left <\,f,\, g_{\,i} \,\right >_{\alpha}\,\right|^{\,2}\,\right\},
\end{align*}
for all \,$f \,\in\, U$, then \,$\left\{\,g_{\,i}\,\right\}_{i \,=\, 1}^{\infty}$\, is a fuzzy \,$K$-frame for \,$U$. 
\end{theorem}

\begin{proof}
For each \,$f \,\in\, U$, we have 
\begin{align*}
&\left(\,\sum\limits_{i \,=\, 1}^{\infty}\, \left|\,\left <\,f,\, f_{\,i} \,\right >_{\alpha}\,\right|^{\,2}\,\right)^{1 \,/\, 2} \,=\, \left\|\,\left\{\,\left <\,f,\, f_{\,i}\,\right >_{\alpha}\,\right\}\,\right\|\\
&\leq\, \left\|\,\left\{\,\left <\,f,\, f_{\,i} \,-\, g_{\,i} \,\right >_{\alpha}\,\right\}\,\right\| \,+\, \left\|\,\left\{\,\left <\,f,\, g_{\,i}\,\right >_{\alpha}\,\right\}\,\right\|\\
&\leq\, \left(\,\sqrt{M} \,+\, 1\,\right)\,\left\|\,\left\{\,\left <\,f,\, g_{\,i}\,\right >_{\alpha}\,\right\}\,\right\|\\
&\Rightarrow\,\sum\limits_{i \,=\, 1}^{\infty}\, \left|\,\left <\,f,\, f_{\,i}\,\right >_{\alpha}\,\right|^{\,2} \,\leq\, \left(\,\sqrt{M} \,+\, 1\,\right)^{\,2}\,\sum\limits_{i \,=\, 1}^{\infty}\, \left|\,\left <\,f,\, g_{\,i}\,\right >_{\alpha}\,\right|^{\,2}\\
&\Rightarrow\,\sum\limits_{i \,=\, 1}^{\infty}\, \left|\,\left <\,f,\, g_{\,i}\,\right >_{\alpha}\,\right|^{\,2} \,\geq\, \dfrac{A}{\left(\,\sqrt{M} \,+\, 1\,\right)^{\,2}}\,\left \|\,K^{\,\ast}\,f\,\right \|_{\alpha}^{\,2}.  
\end{align*}
On the other hand, for each \,$f \,\in\, U$, we have
\begin{align*}
&\left(\sum\limits_{i = 1}^{\infty} \left|\left <f, g_{\,i}\right >_{\alpha}\right|^{\,2}\right)^{1 / 2}\leq\left(\sum\limits_{i = 1}^{\infty} \left|\left <f, f_{\,i} - g_{\,i} \right >_{\alpha}\right|^{\,2}\right)^{1 / 2} + \left(\sum\limits_{i = 1}^{\infty}\, \left|\left <f,\, f_{\,i}\right >_{\alpha}\right|^{\,2}\,\right)^{1 / 2}\\
&\leq\,\left(\,\sqrt{M} \,+\, 1\,\right)\,\left(\,\sum\limits_{i \,=\, 1}^{\infty}\, \left|\,\left <\,f,\, f_{\,i}\,\right >_{\alpha}\,\right|^{\,2}\,\right)^{1 \,/\, 2}\\
&\Rightarrow\,\sum\limits_{i \,=\, 1}^{\infty}\, \left|\,\left <\,f,\, g_{\,i}\,\right >_{\alpha}\,\right|^{\,2} \,\leq\, B\,\left(\,\sqrt{M} \,+\, 1\,\right)^{\,2}\,\left \|\,f\,\right\|_{\alpha}^{\,2}.
\end{align*} 
Thus, \,$\left\{\,g_{\,i}\,\right\}_{i \,=\, 1}^{\infty}$\, is a fuzzy \,$K$-frame for \,$U$.\,This completes the proof.
\end{proof}

\begin{corollary}
Let \,$\left\{\,f_{\,i}\,\right\}_{i \,=\, 1}^{\infty}$\, be a fuzzy frame for \,$U$\, with bounds \,$A$\, and \,$B$.\,Then the following statements are equivalent:
\begin{itemize}
\item[$(i)$]\,$\left\{\,g_{\,i}\,\right\}_{i \,=\, 1}^{\infty}$\, is a fuzzy frame for \,$U$.
\item[$(ii)$]There exists a constant \,$M \,>\, 1$\, such that
\begin{align*}
&\sum\limits_{i \,=\, 1}^{\infty}\, \left|\,\left <\,f,\, f_{\,i} \,-\, g_{\,i}\,\right >_{\alpha}\,\right|^{\,2}\leq\, M\, \min\left\{\,\sum\limits_{i \,=\, 1}^{\infty}\, \left|\,\left <\,f,\, f_{\,i}\,\right >_{\alpha}\,\right|^{\,2},\, \sum\limits_{i \,=\, 1}^{\infty}\, \left|\,\left <\,f,\, g_{\,i} \,\right >_{\alpha}\,\right|^{\,2}\,\right\},
\end{align*}
for all \,$f \,\in\, U$.
\end{itemize}
\end{corollary}

\begin{proof}$(i)\,\Rightarrow\, (ii)$
Suppose \,$\left\{\,g_{\,i}\,\right\}_{i \,=\, 1}^{\infty}$\, is a fuzzy frame for \,$U$\, with bounds \,$C$\, and \,$D$.\,Now, for each \,$f \,\in\, U$, we have
\begin{align*}
&\left(\,\sum\limits_{i \,=\, 1}^{\infty}\, \left|\,\left <\,f,\, f_{\,i} \,-\, g_{\,i} \,\right >_{\alpha}\,\right|^{\,2}\,\right)^{1 \,/\, 2}\\
&\leq\,\left(\,\sum\limits_{i \,=\, 1}^{\infty}\, \left|\,\left <\,f,\, f_{\,i}\,\right >_{\alpha}\,\right|^{\,2}\,\right)^{1 \,/\, 2} \,+\, \left(\,\sum\limits_{i \,=\, 1}^{\infty}\, \left|\,\left <\,f,\, g_{\,i}\,\right >_{\alpha}\,\right|^{\,2}\,\right)^{1 \,/\, 2}\\
&\leq\,\left(\,\sum\limits_{i \,=\, 1}^{\infty}\, \left|\,\left <\,f,\, f_{\,i} \,\right >_{\alpha}\,\right|^{\,2}\,\right)^{1 \,/\, 2} \,+\, \sqrt{D}\,\left \|\,f\,\right \|_{\alpha}\\
&\leq\,\left(\,1 \,+\, \dfrac{\sqrt{D}}{\sqrt{A}}\,\right)\,\left(\,\sum\limits_{i \,=\, 1}^{\infty}\, \left|\,\left <\,f,\, f_{\,i}\,\right >_{\alpha}\,\right|^{\,2}\,\right)^{1 \,/\, 2}.
\end{align*}
On the other hand, for each \,$f \,\in\, U$, we have
\begin{align*}
&\left(\,\sum\limits_{i \,=\, 1}^{\infty}\, \left|\,\left <\,f,\, f_{\,i} \,-\, g_{\,i}\,\right >_{\alpha}\,\right|^{\,2}\,\right)^{1 \,/\, 2}\\
&\leq\,\left(\,\sum\limits_{i \,=\, 1}^{\infty}\, \left|\,\left <\,f,\, f_{\,i}\,\right >_{\alpha}\,\right|^{\,2}\,\right)^{1 \,/\, 2} \,+\, \left(\,\sum\limits_{i \,=\, 1}^{\infty}\, \left|\,\left <\,f,\, g_{\,i}\,\right >_{\alpha}\,\right|^{\,2}\,\right)^{1 \,/\, 2}\\
&\leq\,\sqrt{B}\,\left \|\,f\,\right \|_{\alpha} \,+\, \left(\,\sum\limits_{i \,=\, 1}^{\infty}\, \left|\,\left <\,f,\, g_{\,i}\,\right >_{\alpha}\,\right|^{\,2}\,\right)^{1 \,/\, 2}\\
&\leq\,\left(\,1 \,+\, \dfrac{\sqrt{B}}{\sqrt{C}}\,\right)\,\left(\,\sum\limits_{i \,=\, 1}^{\infty}\, \left|\,\left <\,f,\, g_{\,i}\,\right >_{\alpha}\,\right|^{\,2}\,\right)^{1 \,/\, 2}.
\end{align*}
Now, if we take 
\[M \,=\, \max\left\{\,\left(\,1 \,+\, \dfrac{\sqrt{B}}{\sqrt{C}}\,\right)^{\,2},\, \left(\,1 \,+\, \dfrac{\sqrt{D}}{\sqrt{A}}\,\right)^{\,2}\,\right\},\]
then \,$(ii)$ holds.

$(ii)\,\Rightarrow\,(i)$\,\,Proof directly follows from the proof of the Theorem \ref{3.tmn2.201} by putting \,$K \,=\, I_{F}$.
This completes the proof.    
\end{proof}

\section{Conclusion}
In this paper, we have studied the notion of a fuzzy $K$-frame for fuzzy Hilbert space $U$\, and presented some of their characterizations.\,We have established stability condition of fuzzy \,$K$-frame for \,$U$ under some perturbations.\,We think that these result will be helpful for the researchers to develop fuzzy frame theory.

\end{document}